\documentclass[10pt,a4paper,twoside]{article}

\usepackage{amsmath,amsfonts,amssymb,latexsym,wasysym,mathrsfs}
\usepackage{amsthm}
\usepackage{verbatim}
\usepackage[english]{babel}
\usepackage{a4}
\usepackage{pstricks}

\theoremstyle{plain}
\newtheorem{theorem}{Theorem}[section]
\newtheorem{proposition}{Proposition}[section]
\newtheorem{lemma}{Lemma}[section]

\theoremstyle{definition}
\newtheorem{definition}{Definition}[section]
\theoremstyle{remark}
\newtheorem{remark}{Remark}[section]

\title{The role of disorder in the dynamics of critical fluctuations of mean field models}

\author{Francesca Collet\thanks{e-mail: fcollet@ing.uc3m.es}\\[0.2cm]
{\small Departamento de Ciencia e Ingenier\'ia de Materiales e Ingenier\'ia Qu\'imica} \\
{\small Universidad Carlos III de Madrid}\\
{\small Avenida de la Universidad 30; 28911 - Legan\'es (Madrid), Spain} \\[0.2cm]
{\small and} \\[0.2cm] 
Paolo Dai Pra\thanks{e-mail: daipra@math.unipd.it}\\[0.2cm]
{\small Dipartimento di Matematica Pura ed Applicata} \\
{\small Universit\`a degli Studi di Padova}\\
{\small Via Trieste 63; 35121 - Padova, Italy}}
\date{}

\begin{document}

\maketitle

\begin{abstract}
\noindent The purpose of this paper is to analyze how the disorder affects the dynamics of critical fluctuations for two different types of interacting particle system: the Curie-Weiss and Kuramoto model. 
The models under consideration are a collection of spins and rotators respectively. They both are subject to a mean field interaction and embedded in a site-dependent, i.i.d. random environment. As the number of particles goes to infinity their limiting dynamics become deterministic and exhibit phase transition. The main result concern the fluctuations around this deterministic limit at the critical point in the thermodynamic limit. From a qualitative point of view, it indicates that when disorder is added spin and rotator  systems belong to two different classes of universality, which is not the case for the homogeneous models (i.e., without disorder).

\vspace{0.3cm}

\end{abstract}

\tableofcontents

\section{Introduction}\label{s1}

Interacting particle systems with {\em mean field} interaction are characterized by the complete absence of geometry in the space of configurations, in the sense that the strength of the interaction between particles is independent of their mutual position. The advantage of dealing with this kind of models is that they usually are analytically tractable and it is rather simple derive their macroscopic equations.  Even if the mean field hypothesis may seem too simplistic to describe physical systems, where geometry and short-range interactions are involved, mean field models have been recently applied to social sciences and finance, as in \cite{BrDu01, DPRST09, DaPTo09, DGGaLa, FrBa08, LaLi07}. \\
We briefly introduce the general framework  and some of its peculiar features. 
By {\em mean field} stochastic process we mean a family $x^{(N)} = (x^{(N)}(t))_{t \geq 0}$ with the following characteristics:
\begin{itemize}
\item
 $x^{(N)}(t) = \left(x^{(N)}_1 (t),x^{(N)}_2 (t), \ldots, x^{(N)}_N (t)\right)$ is a Markov process with $N$ components, taking values on a given measurable space $(E,{\cal{E}})$;
\item
Consider the {\em empirical measure}
\[
\rho_N (t) := \frac{1}{N} \sum_{k=1}^N \delta_{x^{(N)}_k(t)},
\]
which is  a random probability on $(E,{\cal{E}})$. Then $(\rho_N (t))_{t \geq 0}$ is a measure-valued {\em Markov} process.
\end{itemize}
Although this is by no means a {\em standard} definition of mean field model, it captures the basic features of the specific models we will consider.

\noindent
Let $(F,{\cal{F}})$ be a topological vector space, and $h: E \rightarrow F$ be a measurable function. Objects of the form
\[
\int h d\rho_N(t) = \frac{1}{N} \sum_{k=1}^N h \left(x^{(N)}_k(t)\right)
\]
are called {\em empirical averages}. In the case the flow $(\int h d\rho_N(t))_{t \geq 0}$ is a Markov process, we say $\int h d\rho_N(t)$ is an {\em order parameter}. Note that the empirical measure itself is an order parameter (taking $F=$ set of signed measures on $E$, and $h(x) = \delta_x$). Whenever possible, it is interesting to find finite dimensional order parameters, i.e. order parameters for which $F$ is finite dimensional. \\
One of the nice aspects of mean field models is that, in many interesting cases, one can prove a {\em Law of Large Numbers} (as $N \rightarrow +\infty$) for the order parameters, and characterize the deterministic limit as a solution of an ordinary differential equation. This limit is often called the {\em McKean-Vlasov} limit. In particular, the differential equation describing the limit evolution of the empirical measure, will be referred to as the  {\em McKean-Vlasov equation}. This equation has the form
\[
\frac{d}{dt} q = {\cal{L}} q,
\]
where ${\cal{L}}$ is a nonlinear operator acting on signed measures on $E$ (even though other spaces may be more convenient for the analysis of ${\cal{L}}$).
\\
Our main interest is the study of the fluctuations of the order parameter around its limiting dynamics. 
We can capture different features of these fluctuations depending on whether or not the time is rescaled with $N$. If time is not rescaled and we consider the evolution in a time interval $[0,T]$, with $T$ fixed, a Central Limit Theorem holds  for the order parameter for all regimes; in other words, the fluctuations of the order parameter converge to a Gaussian process, which is the unique solution of a linear diffusion equation. Whenever time is rescaled in such a way $T$ goes to infinity as $N$ does, we may observe different behaviors. To avoid further complications, we assume the Markov process $x^{(N)}(t)$ has a ``nice'' {\em chaotic} initial condition: $x^{(N)}_1 (0),x^{(N)}_2 (0), \ldots, x^{(N)}_N (0)$ are i.i.d. with common law $q_0(dx)$, where $q_0$ is a stationary, locally stable solution of the McKean-Vlasov equation (the system is in {\em local} equilibrium). 
\begin{itemize}
\item
{\em Subcritical regime}. Suppose $q_0$ is the {\em unique} stationary solution of the McKean-Vlasov equation, and it is {\em linearly stable} (i.e. stable for the linearized equation). Then we expect the Central Limit Theorem holds {\em uniformly in time}; in particular, this provides a Central Limit Theorem for the stationary distribution of $x^{(N)}$. Some results in this direction are shown in \cite{For09}.
\item
{\em Supercritical regime}. Suppose the set of stationary, linearly stable solutions of the McKean-Vlasov equation has cardinality greater than $1$. In this case {\em metastability} phenomena occur at a time scale exponentially growing in $N$.
\item
{\em Critical regime}. This is the case in the boundary of the subcritical regime: denoting by $\mathfrak{L}$ the linearization of ${\cal{L}}$ around $q_0$, the spectrum $\mathrm{Spec}(\mathfrak{L})$ of $\mathfrak{L}$ is contained in $\{z \in \mathbb{C}: Re(z) \leq 0\}$, but there are elements on $\mathrm{Spec}(\mathfrak{L})$ with zero real part. Under a suitable {\em time speed-up}, the elements of the corresponding  eigenspaces may exhibit large and, possibly, non-normal fluctuations (see \cite{Daw83,CoEi88}).
\end{itemize}
Of course the three regimes described above do not cover in general all possibilities, since stable periodic orbits or even stranger attractors may arise. Moreover, the same model could be in different regimes depending on the values of some parameters ({\em phase transition}).

The main subject of this paper is the analysis of the dynamics of the critical fluctuations in disordered mean field models. \\
We consider a mean field model and we add a site-dependent, i.i.d. random environment, acting as an inhomogeneity in the structure of the system; we aim at analyzing the effect of the disorder in the dynamics of critical fluctuations, as compared with the homogeneous case. We deal with the Curie-Weiss and the Kuramoto models.  We are not aware of similar results concerning non-equilibrium critical fluctuations in presence of disorder. Static fluctuations for the random Curie-Weiss model have been studied in \cite{MaPe91}. \\

We now give the basic ideas of how the dynamics of critical fluctuations are determined. As we mentioned above, the deterministic limiting dynamics of the order parameter is described by a nonlinear evolution operator ${\cal{L}}$. The linearization of this equation around a stationary solution gives rise to the so called linearized operator $\mathfrak{L}$. This operator is also related to the normal fluctuation of the process. At the critical point this operator has an eigenvalue with zero real part, while all other elements of the spectrum have negative real part. The eigenspace of the eigenvalue with zero real part  will be called {\em critical direction}, and usually happens to have low dimension: critical phenomena involve the empirical averages corresponding to this subspace.
Thus, our analysis follow the following points.
\begin{itemize}
\item Locating the critical direction.
\item Determining the correct space-time scaling for the critical fluctuations. This requires an approximation of the time evolution of the order parameter that goes beyond the normal approximation.
\item Proving that the rescaled fluctuations vanish along non-critical directions. This will be done using the method of ``collapsing processes'' : it was developed by Comets and Eisele in \cite{CoEi88} for a geometric long-range interacting spin system and was previously applied to a homogeneous mean field spin-flip system in \cite{Sar07}. 
\item Determining the limiting dynamics in the critical direction. It will be done using arguments of perturbation theory for Markov processes, which has been treated in \cite{PaStVa77}, and of tightness, applied to a suitable martingale problem. 
\end{itemize} 


From a qualitative point of view, our results indicate that when disorder is added, spin systems and rotators belong to two different classes of universality, which is not the case for homogeneous systems. Roughly speaking, in spin systems the fluctuations produced by the disorder always prevail in the critical regime: these fluctuations evolve in a time scale of order $N^{\frac{1}{4}}$, while the critical slowing down for homogeneous systems is $N^{\frac{1}{2}}$. For rotators, the disorder does not modify the $N^{\frac{1}{2}}$ slowing down. However, as the ``strength'' of the disorder increases, the Kuramoto model undergoes a further phase transition: for sufficiently small disorder, the dynamics of critical fluctuations converge to a nonlinear, ergodic diffusion, as in the homogeneous case; for larger disorder, the limiting diffusion loses ergodicity, and actually explodes in finite time. \\
We finally remark that in \cite{CDS10} we have analyzed the critical fluctuations for a spin system close in spirit to the Curie-Weiss model, although with a less general disorder distribution.

\section{The Random Curie-Weiss Model}\label{s2}

\subsection{Description of the Model}\label{ss1}

Let $\mathscr{S}=\{-1,+1\}$ be the spin space, and $\mu$ be an even probability on $\mathbb{R}$. Let also \mbox{$\underline{\eta}=(\eta_j)_{j=1}^N \in \mathbb{R}^N$} be a sequence of independent, identically distributed random variables, defined on some probability space $(\Omega, \mathcal{F}, P)$, and distributed according to  $\mu$. They represent a random, inhomogeneous magnetic field.\\
Given a configuration $\underline{\sigma}=(\sigma_j)_{j=1}^N \in \mathscr{S}^N $ and a realization of the magnetic field $\underline{\eta}$, we  define the Hamiltonian $H_N(\underline{\sigma},\underline{\eta}):  \mathscr{S}^N \times \mathbb{R}^N  \rightarrow  \mathbb{R}$ as
\begin{equation}\label{1CWG}   
  H_N(\underline{\sigma},\underline{\eta})=-\frac{\beta}{2N}\sum_{j,k=1}^N \sigma_j \sigma_k - \beta  \sum_{j=1}^N \eta_j \sigma_j  \,,
\end{equation}
where $\sigma_j$ is the spin value at site $j$, and $\eta_j$ is the local magnetic field associated with the same site. Let $\beta>0$  be the inverse  temperature. 
For given $\underline{\eta}$,~ $\underline{\sigma}(t)=(\sigma_j(t))_{j=1}^N$, with $t \geq 0$, is a $N$-spin system evolving as a continuous time Markov chain on $\mathscr{S}^N$, with infinitesimal generator $L_N$ acting on functions 
$f:\mathscr{S}^N \rightarrow \mathbb{R}$ as follows:
\begin{equation}\label{2CWG}
    L_Nf(\underline{\sigma})=\sum_{j=1}^{N} e^{-\beta\sigma_j( m^{\underline{\sigma}}_{N} + \eta_j)}\nabla_j^\sigma f(\underline{\sigma}),
\end{equation}
where $\nabla_j^\sigma f(\underline{\sigma})=f(\underline{\sigma}^j)-f(\underline{\sigma})$ and the $k$-th component of $\underline{\sigma}^j$, which is the  spin flip at the site $j$, is
\begin{displaymath}
    \sigma^j_k=\left\{
    \begin{array}{rcc}
        \sigma_k & \mathrm{for} &k\neq j\\
        -\sigma_k& \mathrm{for} &k= j
    \end{array}
    \right..
\end{displaymath}
The quantity $e^{-\beta\sigma_j( m^{\underline{\sigma}}_{N}+ \eta_j)}$ represents the jump rate of the spins, i.e. the rate at which the transition $\sigma_j\rightarrow -\sigma_j$ occurs for some $j$. The expressions \eqref{1CWG} and \eqref{2CWG} describe a system of mean field ferromagnetically coupled spins, each with its own random magnetic field and subject to Glauber dynamics. The two terms in the Hamiltonian have different effects: the first one tends to align the spins, while the second one tends to point each of them in the direction of its local field.

\begin{remark}
For every value of $\underline{\eta}$, \eqref{2CWG} has a reversible stationary distribution proportional to $\exp[-H_N(\underline{\sigma},\underline{\eta})]$.
\end{remark}

For simplicity, the initial condition $\underline{\sigma}(0)$ is such that $(\sigma_j(0),\eta_j)_{j=1}^N$ are independent and identically distributed with law $\lambda$. Note that, since the marginal law of the $\eta_j$'s is $\mu$, $\lambda$ must be of the form
\begin{equation}\label{initial}
\lambda(\sigma,d\eta) = q_0(\sigma,\eta)\mu(d\eta)
\end{equation}
with $q_0(1,\eta) + q_0(-1,\eta) = 1$, $\mu$-almost surely.
The quantity $\left(\sigma_j(t)\right)_{t \in [0,T]}$ represents the time evolution on $[0,T]$ of $j$-th spin value; it is the trajectory of the single $j$-th spin in time. The space of all these paths is $\mathcal{D}[0,T]$, which is the space of the right-continuous, piecewise-constant functions from $[0,T]$ to $\mathscr{S}$. We endow $\mathcal{D}[0,T]$ with the Skorohod topology, which provides a metric and a Borel $\sigma$-field (see ~\cite{EtKu86} for details).

\subsection{Limiting Dynamics}\label{ss2}

We now describe the dynamics of the process \eqref{2CWG}, in the limit as $N\rightarrow +\infty$, in a fixed time interval $[0,T]$. Later, the equilibrium of the limiting dynamics will be studied. These results are special cases of what shown in  \cite{DaPdHo95}, so proofs are omitted. More details can also be found in  \cite{Col09}.

Let $(\sigma_j [0,T])_{j=1}^N \in (\mathcal{D}[0,T])^N$ denote a path of the system in the time interval $[0,T]$, with $T$ positive and fixed. If $f: \mathscr{S} \times \mathbb{R} \rightarrow \mathbb{R}$, 
we are interested in the asymptotic (as $N \rightarrow +\infty$) behavior of \emph{empirical averages} of the form
\[ 
\frac{1}{N} \sum_{j=1}^{N} f(\sigma_j(t), \eta_j) =: \int f d\rho_N (t)\,,
\]
where $(\rho_N (t))_{t \in [0,T]}$ is the flow of \emph{empirical measures}
\[ 
\rho_N (t) := \frac{1}{N} \sum_{j=1}^{N} \delta_{(\sigma_j (t), \eta_j)} \,.
\]
We may think of $\rho_N:= (\rho_N(t))_{t \in [0,T]}$ as a {\em cadlag} function taking values in $\mathcal{M}_1(\mathscr{S} \times \mathbb{R})$, the space of probability measures on $\mathscr{S} \times \mathbb{R}$ endowed with the weak convergence topology, and the related Prokhorov metric, that we denote by $d_P(\, \cdot \, , \, \cdot \, )$.

The first result we state concerns the dynamics of the flow of empirical measures. We need some more notations. For a given $q: \mathscr{S} \times \mathbb{R} \rightarrow \mathbb{R}$, we introduce the linear operator $\mathcal{L}_q$, acting on $f:\mathscr{S} \times \mathbb{R} \rightarrow \mathbb{R}$ as follows:
\[
\mathcal{L}_q f(\sigma,\eta) := \nabla^{\sigma}\left[ e^{- \beta \sigma \left( m_q + \eta\right)} f(\sigma,\eta) \right],
\]
where
\[
m_q := \int \left[ q(1,\eta) - q(-1,\eta) \right] \mu(d\eta).
\]
Given $\underline{\eta} \in \mathbb{R}^N$, we denote by $\mathcal{P}_N^{\underline{\eta}}$ the distribution on $(\mathcal{D}[0,T])^N$ of the Markov process with generator \eqref{2CWG} and initial distribution $\lambda$. We also denote by 
\[
\mathcal{P}_N\left(d \underline{\sigma}[0,T],d\underline{\eta} \right) := \mathcal{P}_N^{\underline{\eta}}\left(d \underline{\sigma}[0,T]\right)  \mu^{\otimes N}\left(d\underline{\eta} \right)
\]
the joint law of the process and the field.

\begin{theorem}\label{8CWG}
The nonlinear McKean-Vlasov equation
\begin{equation}\label{9CWG}
    \left\{
    \begin{array}{cccr}
        \frac{\partial q_t (\sigma, \eta )}{\partial t} & = & \mathcal{L}_{q_t} q_t (\sigma, \eta) \\
        q_0 (\sigma,\eta) &  & \mbox{given in \eqref{initial}} & \\
    \end{array}
    \right.
\end{equation}
admits a unique solution in $\mathcal{C}^1 \left[[0,T], \left(L^1(\mu) \right)^{\mathscr{S}}\right]$, and $q_t(\cdot,\eta)$ is probability on $\mathscr{S}$, for $\mu$-almost every $\eta$ and every $t>0$. Moreover, for every $\varepsilon >0$ there exists $C(\varepsilon) >0 $ such that
\[
\mathcal{P}_N \left( \sup_{t \in [0,T]} d_P ( \rho_N(t), q_t) > \varepsilon \right) \leq e^{-C(\varepsilon) N}
\]
for $N$ sufficiently large, where, by abuse of notations, we identify $q_t$ with the probability $q_t(\sigma,\eta) \mu(d\eta)$ on $\mathscr{S} \times \mathbb{R}$.
\end{theorem}
Thus, equation \eqref{9CWG} describes the infinite-volume dynamics of the system. Next result gives a characterization of stationary solutions of \eqref{9CWG}.
\begin{lemma}\label{11CWG}
Let $q_*: \mathscr{S} \times \mathbb{R} \rightarrow \mathbb{R}$, such that $q_* (\sigma,\cdot)$ is measurable and $q_*(\cdot,\eta)$ is a probability on $\mathscr{S}$. Then $q_*$ is a stationary solution of \eqref{9CWG}, i.e. $\mathcal{L}_{q_*} q_* \equiv 0$, if and only if it is of the form
\begin{equation}\label{12CWG}
q_*(\sigma,\eta) = \frac{e^{\beta \sigma \left( m_{*} + \eta \right)}}{2 \cosh \left(\beta\left( m_{*} + \eta \right) \right)} \,,
\end{equation}
where $m_*$ satisfies the self-consistency relation
\begin{equation}\label{13CWG}
m_* = \int \left[ q_*(1,\eta) - q_*(-1,\eta) \right] \mu(d\eta).
\end{equation} 
Moreover, $m_*=0$ is always a solution of \eqref{13CWG} and it is linearly (resp. neutrally) stable if and only if 
\begin{equation}\label{linstabcond} 
\beta \int \frac{\mu(d\eta)}{\cosh^2 (\beta \eta)} < \, (\mbox{resp. $=$}) \; 1.
\end{equation}
\end{lemma}

\begin{remark} 
The transition between uniqueness and non-uniqueness of the solution of \eqref{13CWG} in general is not related to the change of stability for $m_*=0$. If the distribution $\mu$ is unimodal on $\mathbb{R}$, the two thresholds coincide: the paramagnetic solution is linearly stable when it is unique and unstable when it is not. In case  we choose $\mu = \frac{1}{2} \left( \delta_{\eta} + \delta_{-\eta} \right)$, with $\eta>0$, the phase diagram is more complex: when \eqref{linstabcond} fails, the paramagnetic solution of \eqref{13CWG} is either unstable, and it coexists with a pair of opposite stable ferromagnetic solutions, or may recover linear stability, coexisting with a pair of unstable ferromagnetic solutions and a pair of stable ferromagnetic ones (see \cite{DaPdHo95} for details). A more general $\mu$ may give rise to arbitrarily many solutions of \eqref{13CWG}.
\end{remark}

\subsection[Dynamics of Critical Fluctuations]{Dynamics of Critical Fluctuations $\left(\beta \displaystyle{\int} \frac{\mu(d\eta)}{\cosh^2 (\beta \eta)} = 1\right)$}\label{ss3}

The results of this section are concerned with the {\em fluctuation flow}
\begin{equation}
\label{fluctflow}
\hat{\rho}_N (t) := \sqrt{N} \left[ \rho_N(t) - q_t \right],
\end{equation}
that takes values on the space of signed measures on $\mathscr{S} \times \mathbb{R}$. It is very convenient to assume that the process starts in {\em local equilibrium}, i.e. $q_0(\sigma,\eta) = q_*(\sigma,\eta)$, where $q_*(\sigma,\eta)$ is a stationary solution of \eqref{9CWG}; it should be not hard to extend all next results to a general initial condition. The proofs of all results stated here will be given in Section \ref{proofsCW}. We first state results valid for all temperatures; later, Lemma \ref{16CWG2}, Proposition \ref{CWpropcritfluct}, Theorems \ref{CWhom} and \ref{CWinhom} are restricted to the critical case.

Functions from $\mathscr{S} \times \mathbb{R}$ are all of the form $F(\sigma, \eta) = \gamma(\eta) + \sigma \phi(\eta)$. However
\[
\int \gamma(\eta) d\hat{\rho}_N (t)  = \sqrt{N} \left[ \frac{1}{N} \sum_{j=1}^N \gamma(\eta_j) - \int \gamma(\eta) \mu(d\eta)\right]
\]
does not change in time, and has a Gaussian limit for every $\gamma \in L^2(\mu)$. Thus, we are only interested in the evolution of integrals of the type
\[
\int \sigma\phi(\eta) d\hat{\rho}_N (t).
\]
It is therefore natural to control the action of the generator $L_N$ on functions of $\underline{\sigma}$ and $\underline{\eta}$ of the form $\psi \left(\int \sigma\phi(\eta) d\hat{\rho}_N \right)$, with 
\[
\hat{\rho}_N :=  \sqrt{N} \left[ \frac{1}{N} \sum_{j=1}^{N} \delta_{(\sigma_j, \eta_j)} - q_* \right].
\]
\begin{proposition}
\label{CWpropfluct}
Let $\psi:\mathbb{R}^n \rightarrow \mathbb{R}$ be of class ${\cal{C}}^1$, and $\phi \in \left(L^2(\nu)\right)^n$, where $\nu$ is the measure on $\mathbb{R}$ defined by
\begin{equation}\label{CWnu}
\nu(d\eta) = \frac{\mu(d\eta)}{\cosh(\beta(m_* + \eta))} . 
\end{equation}
Then
\begin{multline}\label{CWfluctuations}
L_N \psi \left(\int \sigma\phi(\eta) d\hat{\rho}_N \right) \\ = 2\sum_{i=1}^n \partial_i \psi \left(\int \sigma\phi(\eta) d\hat{\rho}_N \right) \left[ \int \sinh(\beta (m_* + \eta)) \phi_i(\eta)  d\hat{\rho}_N - \int \sigma \mathfrak{L}\phi_i(\eta) d\hat{\rho}_N \right] \\ + 2 \sum_{i,j=1}^n \partial^2_{ij} \psi\left(\int \sigma\phi(\eta) d\hat{\rho}_N \right) \int \frac{\phi_i(\eta)\phi_j(\eta)}{\cosh(\beta (m_* + \eta))} \mu(d\eta) + o(1),
\end{multline}
where 
\begin{equation}\label{CWlin}
\mathfrak{L} \phi_i(\eta) = \cosh (\beta (m_* + \eta)) \phi_i(\eta) - \beta \int \frac{\phi_i(\eta)}{\cosh (\beta (m_* + \eta))} \mu(d\eta) \,.
\end{equation}
Moreover the remainder $o(1)$ in (\ref{CWfluctuations}) is of the form 
\begin{equation} \label{CWremform}
R_N\left(\int H(\sigma,\eta) d\hat{\rho}_N \right)
\end{equation}
where $H(\sigma,\eta)$ is the vector-valued function
\begin{multline*}
H(\sigma,\eta) = \left( \sigma\phi(\eta), \sigma, [\cosh(\beta (m_* + \eta)) - \sigma \sinh(\beta (m_* + \eta))]\phi(\eta),\right.  \\ \left. [\sigma \cosh(\beta (m_* + \eta)) - \sinh(\beta (m_* + \eta))]\phi(\eta) \right),
\end{multline*}
and
\begin{equation}\label{CWrem}
\lim_{N \rightarrow +\infty} \sup_{|x|,|y|,|z|,|w| \leq M} R_N(x,y,z,w) = 0
\end{equation}
for every $M>0$.
\end{proposition}

Proposition \ref{CWpropfluct}, whose proof consists of a rather standard computation that will be sketched in Section \ref{proofsCW}, is the essential ingredient for proving a Central Limit Theorem for the empirical flow, i.e. to show that the fluctuation flow converges in law  to a Gaussian process. The proof of this result requires to identify an appropriate Hilbert space for the fluctuations $\hat{\rho}_N$ (see e.g. \cite{CoEi88} for related results). Our main aim is, however, to describe large-time fluctuations at the critical points; the additional technical difficulties arising, have not allowed us to obtained the desired results under the present assumptions, in particular with no requirements on the field distribution $\mu$ (except evenness). Thus we find it preferable to make the following assumption at this point.
\begin{itemize}
\item[{\bf (F)}]
$\mu$ has {\em finite} support $\mathscr{D}$. 
\end{itemize}
Under assumption {\bf (F)}, the space $L^2(\nu)$ is finite-dimensional. Together with the following simple result, this greatly simplifies the analysis of fluctuations.
\begin{lemma}\label{16CWG} 
The operator $\mathfrak{L}$ defined in \eqref{CWlin} is self-adjoint in $L^2(\nu)$. 
\end{lemma} 
Now, for $m := |\mbox{supp}(\mu)|$, let $\varphi_0,\varphi_1,\ldots,\varphi_{m-1}$ be a complete set of eigenvectors for $\mathfrak{L}$, with eigenvalues $\lambda_0\leq \lambda_1 \leq \ldots \leq \lambda_{m-1}$.  Proposition \ref{CWpropfluct}, together with the classical Corollary 8.7, in Chapter 4 of
\cite{EtKu86}, yields the following Central Limit Theorem, whose standard proof is omitted.  
\begin{proposition}\label{CWCLT} 
Set $X^{(N)}_i(t) := \int \sigma \varphi_i (\eta) d \hat{\rho}_N(t)$. 
Then, under $\mathcal{P}_N$, $\left( X^{(N)}_i \right)_{i=0}^{m-1}$ converges in law to the Gaussian process $(X_i)_{i=0}^{m-1}$ solving the following linear stochastic differential equations 
\[ 
dX_i(t) = \left[ \mathscr{H}_i - \lambda_i X_i(t) \right] dt + b_i \, dW_i(t)
\] 
where

\vspace{.2cm}
$\bullet$ 
$\left(X_0(0),X_1(0),\ldots,X_{m-1}(0), \mathscr{H}_0, \mathscr{H}_1,\ldots, \mathscr{H}_{m-1} \right)$ is a centered Gaussian vector with 
{\footnotesize{
\begin{multline*}
\mathrm{Cov} (X_i(0),X_j(0)) = \int \varphi_i(\eta) \varphi_j(\eta) \mu(d\eta) \\
-  \int \varphi_i (\eta)\tanh(\beta (m_*+\eta)) \mu(d\eta) \int \varphi_j (\eta)\tanh(\beta (m_* + \eta)) \mu(d\eta) 
\end{multline*}
\begin{multline*}
\mathrm{Cov}(\mathscr{H}_i,\mathscr{H}_j) = \int \varphi_i (\eta) \varphi_j (\eta) \sinh^2(\beta(m_* + \eta)) \mu(d\eta) \\
-  \int \varphi_i (\eta)\sinh(\beta (m_* + \eta)) \mu(d\eta) \int \varphi_j (\eta)\sinh(\beta (m_* + \eta)) \mu(d\eta)
\end{multline*}
\begin{multline*}
\mathrm{Cov}(\mathscr{H}_i,X_j(0)) = \int \varphi_i (\eta) \varphi_j (\eta) \sinh(\beta (m_* + \eta))\tanh(\beta(m_* + \eta)) \mu(d\eta) \\
-  \int \varphi_i (\eta) \sinh(\beta (m_* + \eta)) \mu(d\eta) \int \varphi_j (\eta)\tanh(\beta (m_* + \eta)) \mu(d\eta)  
\end{multline*}}}

$\bullet$ 
$b_i^2 := \int \varphi_i^2(\eta) \nu(d\eta)$. 

\vspace{.2cm}
$\bullet$ 
$(W_i)_{i=0}^{m-1}$ are independent standard Brownian motions, that are independent of the vector $\left(X_0(0),X_1(0),\ldots,X_{m-1}(0), \mathscr{H}_0, \mathscr{H}_1,\ldots, \mathscr{H}_{m-1} \right)$.

\end{proposition}

Note that the randomness of the field persists in the limiting dynamics of fluctuations, due to the correlated, constant random drifts $\mathscr{H}_i$. Observe that $\mathscr{H}_i \equiv 0 $ if $\mu = \delta_0$, i.e. when the random field is absent. 

We now look more closely at fluctuations around the paramagnetic solution $m_* = 0$ at the {\em critical regime}, i.e. for those values of $\beta$ for which $\beta \int_{\mathscr{D}} \frac{\mu(d\eta)}{\cosh^2(\beta \eta)} = 1$.
\begin{lemma}\label{16CWG2} 
Assume $\beta \int_{\mathscr{D}} \frac{\mu(d\eta)}{\cosh^2(\beta \eta)} = 1$ and $m_* = 0$. Then $\mathfrak{L}$ is nonnegative, and its kernel is spanned by the function $\frac{1}{\cosh(\beta \eta)}$. 
\end{lemma}

In the critical regime $\beta \int_{\mathscr{D}} \frac{\mu(d\eta)}{\cosh^2(\beta \eta)} = 1$, we have $\lambda_0 = 0$, and $\lambda_i > 0$  for $i>0$ (it is actually easily shown that $\lambda_i \geq 1$ for $i>0$). It follows that the process $X_0(t)$ in Proposition \ref{CWCLT} has a variance that diverges as $t \rightarrow +\infty$. A sharper description of the large time fluctuations is obtained by considering more ``moderate" fluctuations:
\[
\tilde{\rho}_N := N^{-\frac{1}{4}} \hat{\rho}_N.
\]
The following result improves the expansion given in Proposition \ref{CWpropfluct}.
\begin{proposition}\label{CWpropcritfluct}
Under the same assumptions of Proposition \ref{CWpropfluct}, and the further conditions $\beta \int_{\mathscr{D}} \frac{\mu(d\eta)}{\cosh^2(\beta \eta)} = 1$ and $m_* = 0$, we have
\begin{multline}
\label{CWcritfluctuations}
L_N \psi \left(\int \sigma\phi(\eta) d\tilde{\rho}_N \right)   =  L^{(1)}\psi \\
+ 2N^{-\frac{1}{4}}\sum_{i=1}^n \partial_i \psi \left(\int \sigma\phi(\eta) d\tilde{\rho}_N \right)   \int \sinh(\beta \eta) \phi_i(\eta)  d\hat{\rho}_N   \\ + N^{-\frac{1}{4}} L^{(2)}\psi + N^{-\frac{1}{2}} L^{(3)}\psi  + o \left( N^{-\frac{1}{2}} \right),
\end{multline}
where
\begin{align*}
& L^{(1)}\psi  := - 2\sum_{i=1}^n \partial_i \psi \left(\int \sigma\phi(\eta) d\tilde{\rho}_N \right) \int \sigma \mathfrak{L}\phi_i(\eta) d\tilde{\rho}_N \\
& L^{(2)}\psi  := -2 \beta \sum_{i=1}^n \partial_i \psi \left(\int \sigma\phi(\eta) d\tilde{\rho}_N \right) \int \sigma d\tilde{\rho}_N  \int \sigma \sinh(\beta \eta) \phi_i(\eta) d\tilde{\rho}_N \\
&L^{(3)}\psi :=  \sum_{i=1}^n \partial_i \psi \left(\int \sigma\phi(\eta) d\tilde{\rho}_N \right) \left[ 2 \beta \int \cosh(\beta \eta) \phi_i(\eta) d\hat{\rho}_N \int \sigma d\tilde{\rho}_N \right. \\ 
& \qquad \left. - \beta^2 \left( \int \sigma d\tilde{\rho}_N \right)^2  \int \sigma \cosh(\beta \eta) \phi_i(\eta) d\tilde{\rho}_N + \frac{\beta^3}{3} \int \frac{\phi_i(\eta)}{\cosh(\beta \eta)} \mu(d\eta) \left( \int \sigma d\tilde{\rho}_N \right)^3 \right] \\ 
& \qquad + 2 \sum_{i,j=1}^n \partial^2_{ij} \psi\left(\int \sigma\phi(\eta) d\hat{\rho}_N \right) \int \frac{\phi_i(\eta)\phi_j(\eta)}{\cosh(\beta \eta)} \mu(d\eta) 
\end{align*}
Moreover the remainder $o \left( N^{-\frac{1}{2}} \right)$ in (\ref{CWcritfluctuations}) is of the form $N^{-\frac{1}{2}} R_N$ with $R_N$ satisfying (\ref{CWrem}).
\end{proposition}

Note that in Proposition \ref{CWpropcritfluct}, functions depending only on $\eta$ are still integrated with respect to $\hat{\rho}$, rather than $\tilde{\rho}$; indeed, by the standard Central Limit Theorem, those integrals with respect to $\hat{\rho}$ have a Gaussian limit under $\mathcal{P}_N$.

Proposition  \ref{CWpropcritfluct} allows to deal easily with the homogeneous case $\mu = \delta_0$. Using the notations of Proposition \ref{CWCLT} we have $m=1$, $\varphi_0 \equiv 1$. Thus, using Proposition \ref{CWpropcritfluct} with $n=1$, $\phi \equiv 1$ and $\beta = \beta_c =1$, we easily observe that $L^{(1)}\psi = L^{(2)}\psi \equiv 0$, and 
\[
L^{(3)}\psi = -\frac{2}{3} \left(\int \sigma d\tilde{\rho}_N\right)^3 \psi'\left(\int \sigma d\tilde{\rho}_N \right) + 2 \psi''\left(\int \sigma d\tilde{\rho}_N \right).
\]
Using convergence of generators as in Proposition \ref{CWCLT} we readily obtain the dynamics of large-time critical fluctuations for the homogeneous model. This result is a simple special case of what obtained in \cite{CoEi88}.
\begin{theorem}\label{CWhom}
Assume $\mu = \delta_0$, and $\beta = 1$. The stochastic process
\[
Y_N(t) := \int \sigma d\tilde{\rho}_N(\sqrt{N}t)
\]
converges weakly, under $\mathcal{P}_N$, to the unique solution of the stochastic differential equation
\[ \left\{
\begin{array}{l}
 dY(t)  =   -\frac{2}{3} \, Y^3(t) \, dt + 2 \, dW(t) \\
 \\
 Y(0) = 0
\end{array} \right.
\]
where $W$ is a standard Brownian motion.
\end{theorem}
As we will see (proofs are in Section \ref{proofsCW}), the inhomogeneous case requires more sophisticated arguments.
\begin{definition}\label{defcollapse}
We say that a sequence of stochastic processes $(\xi_n(t))_n$, for $t\in[0,T]$, \emph{collapses to zero} if for every $\varepsilon >0$,
\[
\lim_{n\rightarrow +\infty} P \left( \sup_{t \in [0,T]} \vert \xi_n(t) \vert > \varepsilon \right) =0
\]
\end{definition}
\begin{theorem}\label{CWinhom}
Assume $m_* =0$, $\beta \int_{\mathscr{D}} \frac{\mu(d\eta)}{\cosh^2(\beta \eta)} = 1$, and, for $i=0,1,\ldots,m-1$, let 
\begin{equation}\label{ordparCW}
Y_i^{(N)}(t) := \int \sigma \varphi_i (\eta) d\tilde{\rho}_N (N^{\frac{1}{4}} t),
\end{equation}
where $\varphi_0,\ldots,\varphi_{m-1}$ is the basis introduced in Proposition \ref{CWCLT}. Under $\mathcal{P}_N$ the processes $\left( Y_i^{(N)}(t) \right)_{i=1}^{m-1}$ collapse to zero, while $Y_0^{(N)}(t)$ converges in law to the process
\[
Y_0(t) := 2\mathscr{H} t,
\]
where $\mathscr{H}$ is  a Gaussian random variable, with zero mean and variance $\int_{\mathscr{D}} \tanh^2(\beta \eta) \mu(d\eta)$.
\end{theorem}
Thus, the disorder has a dramatic impact on fluctuations at the critical points: fluctuations arise at a much shorter time scale ($N^{\frac{1}{4}}$ rather that $N^{\frac{1}{2}}$), and have the simple form of a linear function with random slope.

\section{The Random Kuramoto Model}\label{s3}

\subsection{Description of the Model}\label{ss4}

Let $I=[0,2\pi)$ be the one dimensional torus, and $\mu$ be an even probability on $\mathbb{R}$. Let also $\underline{\eta}=(\eta_j)_{j=1}^N \in \mathbb{R}^N$ be a sequence of independent, identically distributed random variables, defined on some probability space $(\Omega, \mathcal{F}, P)$, and distributed according to $\mu$. Given a configuration \mbox{$\underline{x}=(x_j)_{j=1}^N \in I^N$} and a realization of the random environment $\underline{\eta}$, we can define the Hamiltonian $H_N(\underline{x},\underline{\eta}):  I^N \times \mathbb{R}^N  \rightarrow  \mathbb{R}$ as
\begin{equation}\label{1Ki}   
 H_N(\underline{x},\underline{\eta})=-\frac{\theta}{2N}\sum_{j,k=1}^N \cos(x_k-x_j) + \omega \sum_{j=1}^N  \eta_j  x_j \,,
\end{equation}
where $x_j$ is the position of the rotator at site $j$ and $\omega \eta_j$, with $\omega >0$, can be interpreted as its own frequency. Let $\theta$, positive parameter,  be the coupling strength.
For given $\underline{\eta}$, the stochastic process $\underline{x}(t)=(x_j(t))_{j=1}^N$, with $t \geq 0$, is a \mbox{$N$-rotator} system evolving as a  Markov diffusion process on $I^N$, with infinitesimal generator $L_N$ acting on ${\cal{C}}^2$ functions 
$f:I^N \rightarrow \mathbb{R}$ as follows:
\begin{align}\label{2Ki}
L_Nf(\underline{x})&= \frac{1}{2}  \sum_{j=1}^N\frac{\partial^2 f}{\partial x_j^2}(\underline{x})+\sum_{j=1}^N \frac{\partial H_N}{\partial x_j}(\underline{x},\underline{\eta}) \nonumber\\
&=\frac{1}{2}  \sum_{j=1}^N\frac{\partial^2 f}{\partial x_j^2}(\underline{x})+\sum_{j=1}^N \left\{\omega\eta_j + \frac{\theta}{N} \sum_{k=1}^N \sin(x_k - x_j) \right\}\frac{\partial f}{\partial x_j}(\underline{x}) \,.
\end{align}
Consider the complex quantity
\begin{equation}\label{3Ki}
r_{N} e^{i \Psi_N}=\frac{1}{N}\sum_{j=1}^N e^{i x_j} \, ,
\end{equation}
where $0 \leq r_N \leq 1$ measures the phase coherence of the rotators and $\Psi_N$ measures the average phase. We can reformulate the expression of the infinitesimal generator \eqref{2Ki} in terms of \eqref{3Ki}:
\begin{equation}\label{4Ki}
L_Nf(\underline{x})= \frac{1}{2}  \sum_{j=1}^N\frac{\partial^2 f}{\partial x_j^2}(\underline{x})+\sum_{j=1}^N \left\{\omega\eta_j + \theta r_N  \sin(\Psi_N - x_j) \right\}\frac{\partial f}{\partial x_j}(\underline{x}) \, .
\end{equation}
The expressions \eqref{1Ki} and \eqref{4Ki} describe a system of mean field coupled rotators, each with its own frequency and subject to diffusive dynamics. The two terms in the Hamiltonian have different effects: the first one tends to synchronize the rotators, while the second one tends to make each of them rotate at its own frequency.\\ 

For simplicity, the initial condition $\underline{x}(0)$ is such that $(x_j(0),\eta_j)_{j=1}^N$ are independent and identically distributed with law $\lambda$. We assume $\lambda$ is of the form
\begin{equation}\label{Kinitial}
\lambda(dx,d\eta) = q_0(x,\eta) \mu(d\eta)dx
\end{equation}
with $\int_I q_0(x,\eta) \, dx = 1$, $\mu$-almost surely. The quantity $x_j(t)$ represents the time evolution on $[0,T]$ of $j$-th rotator; it is the trajectory of the single $j$-th rotator in time. The space of all these paths is $\mathcal{C}[0,T]$, which is the space of the continuous function from $[0,T]$ to $I$, endowed with the uniform topology.\\


\subsection{Limiting Dynamics}\label{ss5}

We now describe the dynamics of the process \eqref{2Ki}, in the limit as $N\rightarrow +\infty$, in a fixed time interval $[0,T]$. Later, the equilibrium of the limiting dynamics will be studied. These results are special cases of what shown in  \cite{DaPdHo95}, so proofs are omitted. 

Let $(x_j [0,T])_{j=1}^N \in (\mathcal{C}[0,T])^N$ denote a path of the system in the time interval $[0,T]$, with $T$ positive and fixed. If $f: I \times \mathbb{R} \rightarrow \mathbb{R}$, we are interested in the asymptotic (as $N \rightarrow +\infty$) behavior of \emph{empirical averages} of the form
\[ 
\frac{1}{N} \sum_{j=1}^{N} f(x_j(t), \eta_j) =: \int f d\rho_N (t)\,,
\]
where $(\rho_N (t))_{t \in [0,T]}$ is the flow of \emph{empirical measures}
\[ 
\rho_N (t) := \frac{1}{N} \sum_{j=1}^{N} \delta_{(x_j (t), \eta_j)} \,.
\]
We may think of $\rho_N:= (\rho_N(t))_{t \in [0,T]}$ as a continuous function taking values in $\mathcal{M}_1(I \times \mathbb{R})$, the space of probability measures on $I \times \mathbb{R}$ endowed with the weak convergence topology, and the related Prokhorov metric, that we denote by $d_P(\, \cdot \, , \, \cdot \, )$.

The first result we state concerns the dynamics of the flow of empirical measures. We need some more notations. For a given $q: I \times \mathbb{R} \rightarrow \mathbb{R}$, we introduce the linear operator $\mathcal{L}_q$, acting on $f: I \times \mathbb{R} \rightarrow \mathbb{R}$ as follows:
\begin{equation}\label{16Ki}
\mathcal{L}_q f(x,\eta)= \frac{1}{2}\frac{\partial^2 f}{\partial x^2}(x,\eta) - \frac{\partial}{\partial x} \left\{ \left[ \omega\eta + \theta r_{q} \sin(\Psi_{q} - x) \right] f(x,\eta) \right\},
\end{equation}
where
\[
r_{q} \, e^{i\Psi_{q}} := \int_I \int e^{ix} \, q(x,\eta) \, \mu(d\eta) \, dx.
\]
Given $\underline{\eta} \in \mathbb{R}^N$, we denote by $\mathcal{P}_N^{\underline{\eta}}$ the distribution on $(\mathcal{C}[0,T])^N$ of the Markov process with generator \eqref{2Ki} and initial distribution $\lambda$. We also denote by 
\[
\mathcal{P}_N\left(d \underline{x}[0,T],d\underline{\eta} \right) := \mathcal{P}_N^{\underline{\eta}}\left(d \underline{x}[0,T]\right)  \mu^{\otimes N}\left(d\underline{\eta} \right)
\]
the joint law of the process and the environment.

\begin{theorem}\label{8Ki}
The nonlinear McKean-Vlasov equation
\begin{equation}\label{9Ki}
    \left\{
    \begin{array}{cccr}
        \frac{\partial q_t (x, \eta )}{\partial t} & = & \mathcal{L}_{q_t} q_t (x, \eta) \\
        q_0 (x,\eta) &  & \mbox{given in \eqref{Kinitial}} & \\
    \end{array}
    \right.
\end{equation}
admits a unique solution in $\mathcal{C}^1 \left[[0,T], L^1(dx \otimes \mu) \right]$, and $q_t(\cdot,\eta)$ is probability on $I$, for $\mu$-almost every $\eta$ and every $t>0$. Moreover, for every $\varepsilon >0$ there exists $C(\varepsilon) >0 $ such that
\[
\mathcal{P}_N \left( \sup_{t \in [0,T]} d_P ( \rho_N(t), q_t) > \varepsilon \right) \leq e^{-C(\varepsilon) N}
\]
for $N$ sufficiently large, where, by abuse of notations, we identify $q_t$ with the probability $q_t(x,\eta) \mu(d\eta)dx$ on $I \times \mathbb{R}$.
\end{theorem}

Thus, equation \eqref{9Ki} describes the infinite-volume dynamics of the system. Since $\mu$ is symmetric and the operator $\mathcal{L}$ preserves evenness, we can suppose the average phase $\Psi_{q_t} \equiv 0$, without loss of generality. Next result gives a characterization of stationary solutions of \eqref{9Ki}.

\begin{lemma}\label{11Ki}
Let $q_*:I \times \mathbb{R} \rightarrow \mathbb{R}$, such that $q_*(x, \cdot)$ is measurable and $q_*(\cdot,\eta)$ is a probability on $I$. Then $q_*$ is a stationary solution of \eqref{9Ki}, i.e. $\mathcal{L}_{q_*} q_* \equiv 0$, if and only if it is of the form
\begin{multline}\label{12Ki}
q_*(x,\eta)= (Z_*)^{-1}  \cdot e^{2(\omega\eta x + \theta r_* \cos x) } \left[ e^{4\pi \omega\eta} \int_0^{2\pi} e^{-2(\omega\eta x + \theta r_* \cos x)} dx \right.  \\
 \left. + (1 - e^{4\pi \omega\eta}) \int_0^x e^{-2(\omega\eta y + \theta r_* \cos y)} dy\right] \,,
\end{multline}
where $Z_*$ is a normalizing factor and $r_*$ satisfies the self-consistency relation
\begin{equation}\label{13Ki}
r_* = \int_I \int e^{ix} \, q_*(x,\eta) \, \mu(d\eta) \, dx \, .
\end{equation}
Moreover, $r_*=0$ is always a solution of \eqref{13Ki} and, letting 
\begin{equation}\label{Kcrtval}
\theta_c = \left[ \int \frac{\mu (d\eta)}{1 + 4 (\omega \eta)^2} \right]^{-1} \,,
\end{equation}
we have that
\begin{enumerate}
\item \label{item1} if $\mu$ is unimodal on $\mathbb{R}$, then the solution $r_*=0$ is linearly (resp. neutrally) stable if and only if $\theta < \, (\mbox{resp. $=$}) \; \theta_c \,;$
\item \label{item2} if $\mu= \frac{1}{2} (\delta_{1} + \delta_{-1})$, then the solution $r_*=0$ is linearly (resp. neutrally) stable if and only if $\theta < \, (\mbox{resp. $=$}) \; \theta_c \wedge 2 \,.$
\end{enumerate}
\end{lemma}

\begin{remark} 
The transitions uniqueness/non-uniqueness of the solution of \eqref{13Ki} and stability/instability of $r_*=0$ in general do not occur at the same threshold. It does, however, in the case \ref{item1} of the previous Lemma. The phase diagram related to the case \ref{item2} is more complicated. We refer to \cite{DaPdHo95} for further details.
\end{remark}

\begin{remark}
If $r_* = 0$ the stationary solution \eqref{12Ki} reduces to $q_*(x,\eta) := \frac{1}{2\pi}$.
\end{remark}

\subsection[Dynamics of Critical Fluctuations]{Dynamics of Critical Fluctuations $\left( \theta_c = \left[ \displaystyle{\int \frac{\mu (d\eta)}{1 + 4 (\omega \eta)^2}} \right]^{-1} \right)$}\label{ss6}

The results of this section are concerned with the {\em fluctuation flow}
\begin{equation}
\label{Kfluctflow}
\hat{\rho}_N (t) := \sqrt{N} \left[ \rho_N(t) - q_t \right],
\end{equation}
that takes values on the space of signed measures on $I \times \mathbb{R}$. It is very convenient to assume that the process starts in the particular {\em local equilibrium} $q_0(x,\eta) = q_*(x,\eta) = \frac{1}{2\pi}$, which is the stationary solution of \eqref{9Ki} corresponding to $r_*=0$. The proof of the Central Limit Theorem (Proposition \ref{KCLT}) should be not hard also when $q_*(x,\eta)$ is a {\em sincronous} stationary solution of \eqref{9Ki}, i.e. with  $r_* \neq 0$.
The proofs of all results stated here will be given in Section \ref{proofsK}.\\
If $\phi$ is a function from $I \times \mathbb{R}$, we are interested in the evolution of integrals of the type
\[
\int \phi(x,\eta) d\hat{\rho}_N (t) \,.
\]
It is therefore natural to control the action of the generator $L_N$ on functions of $\underline{x}$ and $\eta$ of the form $\psi \left( \int \phi(x,\eta) d\hat{\rho}_N \right),$ with
\[
\hat{\rho}_N := \sqrt{N} \left[ \frac{1}{N} \sum_{j=1}^N \delta_{(x_j,\eta_j)} - q_* \right] = \sqrt{N} \left[ \frac{1}{N} \sum_{j=1}^N \delta_{(x_j,\eta_j)} - \frac{1}{2\pi} \right] .
\]

\begin{proposition}\label{Kpropfluct}
Let $\psi: \mathbb{R}^n \rightarrow \mathbb{R}$ be of class $\mathcal{C}^2$, and $\phi \in (\mathcal{C}^2([0,2\pi) \times \{-1,1\}))^n$ be $2\pi$-periodic in the first argument. Then
\begin{multline*}
L_N \psi \left( \int \phi(x, \eta) d\hat{\rho}_N \right) = \sum_{i=1}^n \partial_i \psi \left( \int \phi(x, \eta) d\hat{\rho}_N \right) \bigg[ \int \mathfrak{L} \phi_i(x, \eta) d\hat{\rho}_N \\
+ \frac{\theta}{N^{\frac{1}{2}}} \int \frac{\partial \phi_i}{\partial x} (x, \eta) \sin(y-x) d\hat{\rho}_N d\hat{\rho}_N \bigg]
\end{multline*}
\vspace{-0.62cm}
\begin{multline}\label{Kfluctuations}
+ \frac{1}{2} \sum_{i,k=1}^{n} \partial_{ik}^2 \psi \left( \int \phi(x, \eta) d\hat{\rho}_N \right) \bigg[  \int \frac{\partial \phi_i}{\partial x} (x,\eta) \frac{\partial \phi_k}{\partial x} (x,\eta) dq_* \\
+ \frac{1}{N^\frac{1}{2}} \int \frac{\partial \phi_i}{\partial x} (x,\eta) \frac{\partial \phi_k}{\partial x} (x,\eta) d\hat{\rho}_N \bigg]
\end{multline}
where the operator
\begin{multline}\label{Klin}
\!\!\! \mathfrak{L}\phi(x,\eta) = \frac{1}{2} \frac{\partial^2 \phi}{\partial x^2}(x,\eta) + \omega\eta \frac{\partial \phi}{\partial x}(x,\eta) + \theta \bigg[ \cos x \int \cos y \, \phi(y,\eta)\, dq_* \\
+ \sin x \int \sin y \, \phi(y,\eta) \, dq_* \bigg]
\end{multline}
is the linearization of $\mathcal{L}$, given by \eqref{16Ki}, around the equilibrium distribution $q_*$.
\end{proposition}

Unlike the proof of Proposition \ref{CWpropfluct}, which requires an expansion of the generator, Proposition \ref{Kpropfluct} follows by the direct application of the generator; its proof is omitted. It provides the key computation for the proof of the Central Limit Theorem (Proposition \ref{KCLT} below). In order to simplify the analysis, we make the following assumption on the distribution of the random environment.
\begin{itemize}
\item[{\bf (H1)}]
$\mu= \frac{1}{2} (\delta_{1} + \delta_{-1})$ 
\end{itemize}

Because of the structure of the system, it is reasonable to focus on functions from $I \times \mathbb{R}$ of the forms $\phi(x,\eta) = \cos(hx)$, $\sin(hx)$, $\eta\cos(hx)$ or $\eta\sin(hx)$, for $h \geq 1$ integer, and thus on the behavior of 
\[
X_h^{(1,N)}(t) := \int \cos(hx) d\hat{\rho}_N(t), \quad  X_h^{(2,N)}(t) := \int \sin(hx) d\hat{\rho}_N(t),
\]
\[
X_h^{(3,N)}(t) := \int \eta \cos(hx) d\hat{\rho}_N(t) \mbox{ and } X_h^{(4,N)}(t) := \int \eta \sin(hx) d\hat{\rho}_N(t) \,.
\]
Proposition \ref{Kpropfluct}, together with the classical Corollary 8.7, in Chapter 4 of
\cite{EtKu86}, yields the following Central Limit Theorem.

\begin{proposition}\label{KCLT} Assume {\bf (H1)} holds. For $r \geq 1$, consider the following space of sequences
\[
\mathcal{H}_{-r} = \bigg\{{\bf x} = \left( x_h^{(1)}, x_h^{(2)}, x_h^{(3)}, x_h^{(4)} \right)_{h \geq 1} : \| {\bf x } \|_{-r} < +\infty \bigg\},
\]
where
\[
 \| {\bf x } \|^2_{-r} :=
 \sum_{h=1}^{+\infty} \frac{1}{(1+h^2)^r} \left[ \left\vert x_h^{(1)} \right\vert^2 + \left\vert x_h^{(2)} \right\vert^2 
+ \left\vert x_h^{(3)} \right\vert^2 + \left\vert x_h^{(4)} \right\vert^2 \right].
\]
Under $\mathcal{P}_N$, on $\mathcal{H}_{-r}$ the process $\left(X_h^{(1,N)}, X_h^{(2,N)}, X_h^{(3,N)}, X_h^{(4,N)} \right)_{h \geq 1}$ converges in law to the Gaussian process $\left(X_h^{(1)}, X_h^{(2)}, X_h^{(3)}, X_h^{(4)} \right)_{h \geq 1}$ solving the following linear stochastic differential equations
\begin{align*}
dX_h^{(1)}(t) &= \left[ \frac{1}{2} \left( \theta \delta_{1h} - h^2 \right) X_h^{(1)}(t) -h\omega X_h^{(4)}(t) \right] dt + \frac{1}{\sqrt{2}} dW_h^{(1)}(t) \\
dX_h^{(2)}(t) &= \left[ \frac{1}{2} \left( \theta \delta_{1h} - h^2 \right) X_h^{(2)}(t) +h\omega X_h^{(3)}(t) \right] dt + \frac{1}{\sqrt{2}} dW_h^{(2)}(t) \\
dX_h^{(3)}(t) &= \left[ - \frac{h^2}{2} X_h^{(3)}(t) -h\omega X_h^{(2)}(t) \right] dt + \frac{1}{\sqrt{2}} dW_h^{(3)}(t) \\
dX_h^{(4)}(t) &= \left[ - \frac{h^2}{2} X_h^{(4)}(t) +h\omega X_h^{(1)}(t) \right] dt + \frac{1}{\sqrt{2}} dW_h^{(4)}(t) 
\end{align*}
where $\delta_{1h}$ is Kronecker delta and

\vspace{.2cm}
$\bullet$ 
$\left(X_h^{(1)}(0), X_h^{(2)}(0), X_h^{(3)}(0), X_h^{(4)}(0) \right)_{h \geq 1}$ is a centered Gaussian vector with \mbox{$\mathrm{Cov} \left( X_h^{(i)}(0), X_k^{(j)}(0) \right) = 0$} for $i \neq j$ or $h \neq k$ and $\mathrm{Var} \left( X_h^{(i)}(0) \right) = \frac{1}{2}$ for any $i,h$.

\vspace{.2cm}
$\bullet$ 
$\left( W_h^{(1)}, W_h^{(2)}, W_h^{(3)}, W_h^{(4)} \right)_{h \geq 1}$ are independent standard Brownian motions, that are independent of $\left(X_h^{(1)}(0), X_h^{(2)}(0), X_h^{(3)}(0), X_h^{(4)}(0) \right)_{h \geq 1}$.
\end{proposition}
Note that the randomness of the field appears only through the parameter $\omega$ in the dynamics of fluctuations. The only source of stochasticity is due to the Brownian motions.

We now proceed to the analysis of the critical regime, i.e. for $\theta = \theta_c \wedge 2$, where $\theta_c$ is given in (\ref{Kcrtval}) and, under {\bf (H1)}, $\theta_c = 1+4 \omega^2$. We make the following further assumption.
\begin{itemize}
\item[{\bf (H2)}]
$\omega < \frac{1}{2}$.
\end{itemize}

Under assumptions {\bf (H1)-(H2)}, we have sufficient control  of the spectrum of $\mathfrak{L}$, as operator in $L^2([0,2\pi) \times \{-1,1\})$. In particular, $\mathfrak{L}$ can be diagonalized in the critical regime, as stated in next Lemma.

\begin{lemma}\label{lmm:specK}
Under assumptions {\bf (H1)-(H2)}, $\ker (\mathfrak{L}) \neq 0$ if and only if $\theta = 1 + 4\omega^2.$ In this last case the spectrum of $\mathfrak{L}$ is given by
\[
\mathrm{Spec}(\mathfrak{L}) = \left\{0, -\frac{1}{2} + 2\omega^2\right\} \cup \left\{-\frac{k^2}{2} \pm i k\omega, \, k \in \mathbb{Z} \setminus \{-1,0,+1\} \right\} \,,
\]
with corresponding eigenspaces 
\[
\begin{array}{rcl}
\ker (\mathfrak{L}) & = & {\mbox span}\left(v_1^{(1)},v_1^{(2)} \right) \\
Eig\left( -\frac{1}{2} + 2\omega^2 \right) & = &  {\mbox span}\left(v_1^{(3)},v_1^{(4)}\right) \\
Eig\left(-\frac{k^2}{2} + i k\omega\right) & = &  {\mbox span}\left(v_k^{(1)},v_k^{(2)}\right) \\
Eig\left(-\frac{k^2}{2} - i k\omega\right) & = &  {\mbox span}\left(v_k^{(3)},v_k^{(4)}\right),
\end{array}
\]
where
\begin{equation}\label{basis}
\begin{array}{lll}
v_1^{(1)}(x,\eta) := \cos x - 2 \omega \eta \sin x && v_1^{(2)}(x,\eta) := \sin x + 2 \omega \eta \cos x \\
 v_1^{(3)}(x,\eta) := \eta \cos x + 2 \omega \sin x && v_1^{(4)}(x,\eta) := 2 \omega \cos x - \eta \sin x \\
 v_k^{(1)}(x,\eta) := \sin (kx) - i \eta \cos (kx) && v_k^{(2)}(x,\eta) := \cos (kx) + i \eta \sin (kx) \\
 v_k^{(3)}(x,\eta) :=  \sin (kx) + i \eta \cos (kx) && v_k^{(4)}(x,\eta) := \cos (kx) - i \eta \sin (kx).
\end{array}
\end{equation}
\end{lemma}

In the critical regime $\theta = \theta_c = 1 + 4\omega^2$ the variance of the processes
\[
U^{(1,N)}(t) := X_1^{(1,N)}(t) - 2 \omega X_1^{(4,N)}(t)  \mbox{ and } U^{(2,N)}(t) := X_1^{(2,N)}(t) + 2 \omega X_1^{(3,N)}(t) \,,
\]
which are the fluctuations of the empirical averages corresponding to the directions generating the kernel of operator $\mathfrak{L},$ diverge as $t \to +\infty.$ A sharper description of the large time fluctuations is obtained by considering more ``moderate'' fluctuations:
\[
\tilde{\rho}_N := N^{-\frac{1}{4}} \hat{\rho}_N.
\]
We will obtain asymptotics, as $N \rightarrow +\infty$, for the signed measures $\tilde{\rho}_N(\sqrt{N} t)$. Note that these measures are completely characterized by their integrals
\begin{equation} \label{Korderpar}
V_h^{(i,N)} (t) := \int v_h^{(i)}(x,\eta) \, d \tilde{\rho}_N(\sqrt{N} t),
\end{equation}
with $h \geq1$ and $i=1,2,3,4$.


\begin{theorem}\label{Kinhom}
Assume $\theta_c = 1 + 4\omega^2$, and $\omega \leq \frac{1}{2 \sqrt{2}}$.
Under $\mathcal{P}_N$ the processes $\left( V_h^{(i,N)}(t) \right)_{h \geq 2}$, for $i=1,2,3,4,$ and $V_1^{(3,N)}(t), V_1^{(4,N)}(t)$ collapse to zero in the sense of Definition \ref{defcollapse}, while the process $\left( V_1^{(1,N)}(t), V_1^{(2,N)}(t) \right)$ converges weakly to the unique solution $\left(V^{(1)}(t), V^{(2)}(t)\right)$ of the stochastic differential equation
\[ \left\{
\begin{array}{l}
 dV^{(1)}(t) =  - \frac{(1+4\omega^2)^2(1-8\omega^2)}{4(1-4\omega^2)^3(1+\omega^2)} \, V^{(1)}(t) \, \left[\left( V^{(1)}(t) \right)^2 + \left( V^{(2)}(t) \right)^2 \right] \, dt + \sqrt{\frac{1+4\omega^2}{2}} \, dW^{(1)}(t)\\
\\
 dV^{(2)}(t)  =  - \frac{(1+4\omega^2)^2(1-8\omega^2)}{4(1-4\omega^2)^3(1+\omega^2)} \, V^{(2)}(t) \, \left[\left( V^{(1)}(t) \right)^2 + \left( V^{(2)}(t) \right)^2 \right]\, dt + \sqrt{\frac{1+4\omega^2}{2}} \, dW^{(2)}(t)\\
\\
 V^{(1)}(0) = V^{(2)}(0) = 0 
\end{array} \right.
\]
where $W^{(1)}$ and $W^{(2)}$ are two independent standard Brownian motions. \\ In the case $ \frac{1}{2 \sqrt{2}} < \omega < \frac{1}{2}$, the process $\left(V^{(1)}(t), V^{(2)}(t)\right)$ explodes in finite time; the convergence above holds for the {\em localized} processes: for every $r>0$, the process $\left( V_1^{(1,N)}(t \wedge T_{N,r}), V_1^{(2,N)}(t \wedge T_{N,r}) \right)$ converges weakly to $\left(V^{(1)}(t \wedge T_r), V^{(2)}(t \wedge T_r)\right)$, where
\[
T_{N,r} := \inf\left\{t>0: \left( V_1^{(1,N)}(t)\right)^2 + \left(V_1^{(2,N)}(t)\right)^2 \geq r \right\} 
\]
\[
 T_r := \inf\left\{t>0: \left( V^{(1)}(t) \right)^2 + \left( V^{(2)}(t) \right)^2 \geq r \right\} .
\]
\end{theorem}

By Theorem \ref{Kinhom} we can derive the limiting dynamics of the critical fluctuations for the homogeneous model $\mu = \delta_0.$ They can be obtained as a particular case setting $\omega=0$.

\begin{theorem}\label{Khom}
Assume $\theta_c = 1.$ For $h \geq 1$ integer, let
\[
Y_h^{(1,N)}(t) := \int \cos(hx) d\tilde{\rho}_N(\sqrt{N}t) \mbox{ and }  Y_h^{(2,N)}(t) := \int \sin(hx) d\tilde{\rho}_N(\sqrt{N}t) \,.
\]
Under $\mathcal{P}_N$ the processes $\left( Y_h^{(i,N)}(t) \right)_{h \geq 2}$, for $i=1,2,$ collapse to zero in the sense of Definition \ref{defcollapse}, while the process $\left( Y_1^{(1,N)}(t), Y_1^{(2,N)}(t) \right)$ converges weakly to the unique solution of the stochastic differential equation
\[ \left\{
\begin{array}{l}
 dY^{(1)}(t)  = -\frac{1}{4} \, Y^{(1)}(t) \, \left[\left( Y^{(1)}(t) \right)^2 + \left( Y^{(2)}(t) \right)^2 \right] \, dt + \frac{1}{\sqrt{2}} \, dW^{(1)}(t)\\
\\
 dY^{(2)}(t)  =   -\frac{1}{4} \, Y^{(2)}(t) \, \left[ \left( Y^{(1)}(t) \right)^2 + \left( Y^{(2)}(t) \right)^2 \right]\, dt + \frac{1}{\sqrt{2}} \, dW^{(2)}(t)\\
\\
Y_1^{(1)}(0)  = Y_1^{(2)}(0)  =0 
\end{array} \right.
\]
where $W^{(1)}$ and $W^{(2)}$ are two independent standard Brownian motions.
\end{theorem}


\section{Collapsing processes}\label{collapsing}

Before giving the details of the proofs of the results stated previously, we briefly present one of the key technical tool: a Lyapunov-like condition, that guarantees a rather strong form of convergence to zero of a sequence of stochastic processes. The first result (Proposition \ref{3}) we state concerns semimartingales driven by Poisson processes, whose proof can be found in the Appendix of \cite{CoEi88}. In the case where the driving noises are Brownian motions, the result takes a slightly simpler form (Proposition \ref{3cont}); its proof is a simple adaptation of the one in \cite{CoEi88}, and it is omitted.

\begin{proposition}\label{3}
Let $\{\xi_n(t)\}_{n \geq 1}$ be a sequence of positive semimartingales on a probability space $(\Omega, \mathscr{A}, \mathscr{P})$, with
\[
d\xi_n(t) = S_n(t)dt + \int_{\mathscr{Y}} f_n(t^-, y) [\Lambda_n(dt,dy) - A_n(t,dy)dt].
\]
Here, $\Lambda_n$ is a Point Process of intensity $A_n(t,dy)dt$ on $\mathbb{R}^+\times \mathscr{Y}$, where $\mathscr{Y}$ is a measurable space, and $S_n(t)$ and $f_n(t)$ are $\mathscr{A}_t$-adapted processes, if we consider $(\mathscr{A}_t)_{t\geq 0}$ a filtration on $(\Omega, \mathscr{A}, \mathscr{P})$ generated by $\Lambda_n$.\\ 
Let $d>1$ and $C_i$ constants independent of $n$ and $t$. Suppose $\{\kappa_n\}_{n \geq1}$, $\{\alpha_n\}_{n \geq1}$ and $\{\beta_n\}_{n \geq 1}$, increasing sequences with
\begin{equation} \label{con1}\tag{$a1$}
\kappa_n^{\frac{1}{d}}\alpha_n^{-1} \xrightarrow{n\rightarrow +\infty} 0, \, \kappa_n^{-1} \alpha_n  \xrightarrow{n\rightarrow +\infty} 0, \, \kappa_n^{-1}\beta_n  \xrightarrow{n\rightarrow +\infty} 0 
\end{equation}
and
\begin{equation} \label{con2}\tag{$a2$}
E \Big[ \Big( \xi_n(0) \Big)^d \Big] \leq C_1 \alpha_n^{-d} \qquad \mbox{for all $n$} \,. 
\end{equation}
Furthermore, let $\{\tau_n\}_{n \geq 1}$ be stopping times such that for $t \in [0, \tau_n]$ and $n \geq 1$,
\begin{equation} \label{con3}\tag{$a3$}
S_n(t) \leq -\kappa_n \delta \xi_n(t) + \beta_n C_2 + C_3 \qquad \mbox{with $\delta>0$,} 
\end{equation}
\begin{equation} \label{con4}\tag{$a4$}
\sup_{\omega\in\Omega, y\in\mathscr{Y}, t\leq\tau_n}\vert f_n(t,y) \vert \leq C_4 \alpha_{n}^{-1} \, ,
\end{equation}
\begin{equation}\label{con5}\tag{$a5$}
\int_{\mathscr{Y}} (f_n(t,y))^2 A_n(t,dy) \leq C_5\, .
\end{equation}
Then, for any $\varepsilon >0$, there exist $C_6>0$ and $n_0$ such that
\begin{equation}\label{3b}
\sup_{n \geq n_0} \mathscr{P} \left\{ \sup_{0 \leq t \leq T \wedge \tau_n} \xi_n (t) > C_6 \left( \kappa_n^{\frac{1}{d}} \alpha_n^{-1} \vee \alpha_n \kappa_n^{-1} \right) \right\} \leq \varepsilon \,.
\end{equation}
\end{proposition}

\begin{proposition}\label{3cont}
Let $\{\xi_n(t)\}_{n \geq 1}$ be a sequence of positive semimartingales on a probability space $(\Omega, \mathscr{A}, \mathscr{P})$, with
\[
d\xi_n(t) = S_n(t)dt + \sum_{i=1}^{m_n} f_n(t,i) dW_i(t)\,.
\]
Here, $\left(W_i \right)_{i=1}^{m_n}$ are independent standard Brownian motions which generate a filtration $\left(\mathscr{A}_t \right)_{t \geq 0}$, and $S_n(t)$ and $f_n(t,i)$ are $\mathscr{A}_t$-adapted processes.\\ 
Let $d>1$ and $C_i$ constants independent of $n$ and $t$. Suppose $\{\kappa_n\}_{n \geq1}$, $\{\alpha_n\}_{n \geq1}$ and $\{\beta_n\}_{n \geq 1}$, increasing sequences with
\begin{equation} \label{con1cont}\tag{$b1$}
\kappa_n^{\frac{1}{d}}\alpha_n^{-1} \xrightarrow{n\rightarrow +\infty} 0, \, \kappa_n^{-1} \alpha_n  \xrightarrow{n\rightarrow +\infty} 0, \, \kappa_n^{-1}\beta_n  \xrightarrow{n\rightarrow +\infty} 0 
\end{equation}
and
\begin{equation} \label{con2cont}\tag{$b2$}
E \Big[ \Big( \xi_n(0) \Big)^d \Big] \leq C_1 \alpha_n^{-d} \qquad \mbox{for all $n$} \,. 
\end{equation}
Furthermore, let $\{\tau_n\}_{n \geq 1}$ be stopping times such that for $t \in [0, \tau_n]$ and $n \geq 1$,
\begin{equation} \label{con3cont}\tag{$b3$}
S_n(t) \leq -\kappa_n \delta \xi_n(t) + \beta_n C_2 + C_3 \qquad \mbox{with $\delta>0$,} 
\end{equation}
\begin{equation}\label{con5cont}\tag{$b4$}
\sum_{i=1}^{m_n} f_n(t,i)^2 \leq C_5\, .
\end{equation}
Then, for any $\varepsilon >0$, there exist $C_6>0$ and $n_0$ such that
\begin{equation}\label{3bcont}
\sup_{n \geq n_0} \mathscr{P} \left\{ \sup_{0 \leq t \leq T \wedge \tau_n} \xi_n (t) > C_6 \left( \kappa_n^{\frac{1}{d}} \alpha_n^{-1} \vee \alpha_n \kappa_n^{-1} \right)\right\} \leq \varepsilon \,.
\end{equation}
\end{proposition}

\section{Proofs for the Random Curie-Weiss Model}\label{proofsCW}

\subsection{Preliminaries}

\paragraph{Proof of Lemma \ref{11CWG}.} An equilibrium probability density for \eqref{9CWG} must satisfy 
\[
\nabla^\sigma \left[ e^{ -\beta \sigma \left( m_{q_t} + \eta \right)} q_t(\sigma,\eta)\right] = 0 \,,
\] 
which is equivalent to $e^{ -\beta \sigma \left( m_{*} + \eta \right)} q_*(\sigma,\eta) = e^{ \beta \sigma \left( m_{*} + \eta \right)} q_*(-\sigma,\eta),$ where $m_*$ is defined in \eqref{13CWG}. Solving, we obtain 
\[
q_*(\sigma,\eta)= e^{ \beta \sigma \left( m_{*} + \eta \right)}\,,
\] 
with the normalizing constant
\[
Z_*=\int_{\mathscr{S}} e^{ \beta \sigma \left( m_{*} + \eta \right)} d\sigma = 2 \cosh \left(\beta\left( m_{*} + \eta \right) \right) \,,
\]
and the proof is complete.

\paragraph{Proof of Lemma \ref{16CWG}.} Obviously $\mathfrak{L}$ is a linear and continuous operator. We have to prove that, if $\phi_1, \phi_2 \in L^2(\nu)$, then $\int_{\mathscr{D}} \left( \mathfrak{L} \phi_1(\eta) \right) \phi_2(\eta) \nu(d\eta) = \int_{\mathscr{D}} \phi_1(\eta) \left( \mathfrak{L} \phi_2(\eta) \right) \nu(d\eta)$. Thus,
\begin{allowdisplaybreaks}
\begin{align*}
\int_{\mathscr{D}} & \left( \mathfrak{L} \phi_1(\eta) \right) \phi_2(\eta) \nu(d\eta) \\
&= \int_{\mathscr{D}} \left[ \cosh(\beta (m_* + \eta)) \phi_1(\eta) - \beta \int_{\mathscr{D}} \frac{\phi_1(\eta)}{\cosh(\beta (m_* + \eta))}\mu(d\eta) \right] \phi_2(\eta) \nu(d\eta)\\
&= \int_{\mathscr{D}} \left[ \cosh(\beta (m_* + \eta)) \phi_2(\eta) - \beta \int_{\mathscr{D}} \frac{\phi_2(\eta)}{\cosh(\beta (m_* + \eta))}\mu(d\eta) \right] \phi_1(\eta) \nu(d\eta)\\
&=\int_{\mathscr{D}} \phi_1(\eta) \left( \mathfrak{L} \phi_2(\eta) \right) \nu(d\eta)
\end{align*}
\end{allowdisplaybreaks}\\
and the proof of self-adjointness is completed. 

\paragraph{Proof of Lemma \ref{16CWG2}.} To prove positivity of $\mathfrak{L}$ we have to show that if $\phi \in L^2(\nu)$, then \mbox{$\int_{\mathscr{D}} \left(\mathfrak{L} \phi(\eta) \right) \phi(\eta) \nu(d\eta) \geq 0$}. Indeed we have
\begin{allowdisplaybreaks}
\begin{align*}
\int_{\mathscr{D}} \left(\mathfrak{L} \phi(\eta) \right) \phi(\eta) \nu(d\eta) &= \int_{\mathscr{D}} \left[ \cosh(\beta \eta) \phi(\eta) - \beta \int_{\mathscr{D}} \frac{\phi(\eta)}{\cosh(\beta \eta)}\mu(d\eta) \right] \phi(\eta) \nu(d\eta)\\
&= \frac{1}{\beta} \int_{\mathscr{D}} \cosh^2(\beta\eta) \phi^2(\eta) \frac{\beta}{\cosh^2(\beta \eta)} \mu(d\eta) \\
& \qquad - \frac{1}{\beta} \left( \int_{\mathscr{D}} \cosh(\beta \eta) \phi(\eta) \frac{\beta}{\cosh^2(\beta \eta)} \mu(d\eta) \right)^2  \geq  0 \,,
\end{align*}
\end{allowdisplaybreaks}\\
where we have used Jensen's inequality for the probability \mbox{$\displaystyle{\beta}$ $ \frac{\mu(d\eta) }{\cosh^2 (\beta \eta)} $}. Moreover, equality holds true if and only if  $\cosh(\beta \eta)\phi(\eta)$ is constant; therefore the null space of the operator $\mathfrak{L}$ is generated by the functions of the form $\phi(\eta)=\frac{1}{\cosh(\beta \eta)}$.

\subsection{Expansions of the Infinitesimal Generator}

\paragraph{Proof of Proposition \ref{CWpropfluct}.} 


By direct computation, and Taylor expansion of $\psi$, we obtain
\begin{allowdisplaybreaks}
\begin{align*}
&L_N\psi \left(\int \sigma \phi (\eta) d\hat{\rho}_N \right)\\
&= \sum_{j=1}^N e^{-\beta \sigma_j \left( \frac{1}{\sqrt{N}} \int \sigma d\hat{\rho}_N + m_* + \eta_j\right)} \left[ \psi \left( \int \sigma \phi (\eta) d\hat{\rho}_N - \frac{2\sigma_j }{\sqrt{N}} \phi(\eta_j) \right) - \psi \left( \int \sigma \phi (\eta) d\hat{\rho}_N \right)\right] \\
&= \sum_{j=1}^N \exp \{ -\beta \sigma_j (m_* + \eta_j)\} \left[ 1 + \sum_{h=1}^3 \frac{1}{h!}\left( -\frac{\beta \sigma_j}{\sqrt{N}} \int \sigma d\hat\rho_N\right)^h + o \left( \frac{1}{N^{\frac{3}{2}}}\right) \right] \\
&\qquad \times \left[ - \frac{2\sigma_j}{\sqrt{N}} \sum_{i=1}^{n} \partial_i \psi(\cdot) \phi_i(\eta_j) + \frac{2}{N} \sum_{i,k=1}^{n} \partial^2_{ik} \psi(\cdot) \phi_i(\eta_j) \phi_k(\eta_j) +o \left( \frac{1}{N}\right)\right] \\
& = - \frac{2}{\sqrt{N}} \sum_{i=1}^n \partial_i \psi(\cdot)  \left\{ \sum_{j=1}^N \phi_i(\eta_j) [\sigma_j \cosh(\beta (m_* + \eta_j)) - \sinh(\beta (m_* + \eta_j))] \right. \\
&\qquad - \frac{\beta}{\sqrt{N}} \left( \int \sigma d\hat{\rho}_N \right) \sum_{j=1}^N  \phi_i(\eta_j) [\cosh(\beta (m_* + \eta_j))  -  \sigma_j \sinh(\beta (m_* + \eta_j))] \\
&\qquad + \frac{\beta^2}{2N} \left( \int \sigma d\hat{\rho}_N \right)^2 \sum_{j=1}^N \phi_i(\eta_j) [\sigma_j \cosh(\beta (m_* + \eta_j)) - \sinh(\beta (m_* + \eta_j))] \\
&\qquad \left. - \frac{\beta^3}{6N^{\frac{3}{2}}} \left( \int \sigma d\hat{\rho}_N \right)^3 \sum_{j=1}^N  \phi_i(\eta_j) [\cosh(\beta (m_* + \eta_j))  -  \sigma_j \sinh(\beta (m_* + \eta_j))]  \right\}\\
& + \frac{2}{N} \sum_{i,k=1}^n \partial_{ik}^2 \psi(\cdot)  \left\{ \sum_{j=1}^N \phi_i(\eta_j) \phi_k(\eta_j) [\cosh(\beta (m_* + \eta_j))  - \sigma_j \sinh(\beta (m_* + \eta_j))] \right\} + o (1) 
\end{align*}
\end{allowdisplaybreaks}\\
We now represent all the terms as integrals with respect to the measure $\hat{\rho}_N$; since 
\[
\int \left[-\sigma\cosh(\beta (m_* + \eta)) + \sinh(\beta (m_* + \eta)) \right] \phi_i(\eta) q_*(d\sigma,d\eta) =0
\] 
and 
\[
\int \left[\cosh(\beta (m_* + \eta)) -\sigma\sinh(\beta (m_* + \eta)) \right] \phi_i(\eta) q_*(d\sigma,d\eta)  = \int_{\mathscr{D}}\frac{\phi(\eta)}{\cosh(\beta (m_* + \eta))} \mu(d\eta)\,,\] 
we obtain
\begin{allowdisplaybreaks}
\begin{align*}
&L_N\psi \left(\int \sigma \phi (\eta) d\hat{\rho}_N \right) \\
&= 2 \sum_{i=1}^n \partial_i \psi( \cdot) \bigg\{  - \int \sigma \left[ \cosh(\beta (m_* + \eta)) \phi_i(\eta) - \beta \int_{\mathscr{D}} \frac{\phi_i(\eta)}{\cosh(\beta (m_* + \eta))} \mu(d\eta)\right] d\hat{\rho}_N \\ 
&\qquad + \int \sinh(\beta (m_* + \eta)) \phi_i(\eta) d\hat{\rho}_N  \\
&\qquad + \frac{\beta}{\sqrt{N}} \int \sigma d\hat{\rho}_N \int [ \cosh(\beta (m_* + \eta)) - \sigma \sinh(\beta (m_* + \eta)) ] \phi_i(\eta) d\hat{\rho}_N \\
&\qquad - \frac{\beta^2}{2N} \! \left( \int \sigma d\hat{\rho}_N \right)^2 \!\!\! \int [\sigma \cosh(\beta (m_* + \eta)) - \sinh(\beta (m_* + \eta))] \phi_i(\eta) d\hat{\rho}_N \\
&\qquad + \frac{\beta^3}{6N}  \left( \int \sigma d\hat{\rho}_N \right)^3 \int_{\mathscr{D}} \frac{\phi_i(\eta)}{\cosh(\beta (m_* + \eta))} \mu(d\eta)  \\
&\qquad + \frac{\beta^3}{6N^{\frac{3}{2}}} \left( \int \sigma d\hat{\rho}_N \right)^3 \int [\cosh(\beta (m_* + \eta)) - \sigma \sinh(\beta (m_* + \eta))] \phi_i(\eta) d\hat{\rho}_N \bigg\} \\
&\quad +2 \sum_{i,k=1}^n \partial_{ik}^2 \psi(\cdot)  \left\{ \int_{\mathscr{D}} \frac{\phi_i(\eta)\phi_k(\eta)}{\cosh(\beta (m_* + \eta))} \mu(d\eta) \right\}  + o (1) \,,
\end{align*}
\end{allowdisplaybreaks}\\
from which \eqref{CWfluctuations} follows. The fact that the remainder $o(1)$ has the form (\ref{CWremform}) and satisfies (\ref{CWrem}) is implied by the Lagrange form of the remainder of the Taylor expansions we have used.

\paragraph{Proof of Proposition \ref{CWpropcritfluct}.} 
It is obtained by a simple rescaling of the last expansion of $L_N \psi\left( \int \sigma \phi(\eta)d\tilde{\rho}_N \right)$ seen in the proof of Proposition \ref{CWpropfluct}. The details are omitted.

\subsection{Collapsing Terms}\label{ssCW:Collapse}

For $N \geq 1$, $M>0$ define the family of stopping times
\[
\tau_{N}^M:=\inf_{t \geq 0} \, \bigg \{ \left \vert Y_i^{(N)}(t) \right \vert \geq M \quad \mbox{for at least a value of $i=0, \dots, m-1$} \bigg \},
\]
where the $Y_i^{(N)}$'s have been defined in (\ref{ordparCW}).
In the rest of this section, we often consider the {\em time-rescaled} infinitesimal generator $J_N=N^{\frac{1}{4}} L_N$, where $L_N$ is given by \eqref{CWcritfluctuations}. 
Whenever we write
\[
J_N \psi\left(Y_0^{(N)}, Y_1^{(N)}, \ldots, Y_{m-1}^{(N)}\right)(t) \,,
\]
we mean
\[
J_N \psi\left(\int \sigma \varphi_0 d\tilde{\rho}_N, \int \sigma \varphi_1 d\tilde{\rho}_N, \dots, \int \sigma \varphi_{m-1} d\tilde{\rho}_N \right) \left\vert_{\tilde{\rho}_N = \tilde{\rho}_N \left(N^{\frac{1}{4}} t\right)}. \right.
\]
We later consider, for $j \in \mathscr{S}$ and $k \in \mathscr{D}$, the counting process $\Lambda^{\sigma}_N (j,k,t)$ which counts the number of spin flips of spins $\sigma_i$ such that $\sigma_i = j$ and $\eta_i = k$, up to time $N^{\frac{1}{4}} t$.
We consider the following semi-martingale decomposition
\begin{equation}\label{CWmartdec}
d\left(Y_i^{(N)}(t)\right)^2 = J_N \left[\left(Y_i^{(N)}\right)^2 \right](t)\,dt+d\mathcal{M}_{N,Y_i^2}^t \,, 
\end{equation}
with $\mathcal{M}_{N,Y_i^2}^t$  the local martingale  given by
\begin{equation} \label{CWmartterm}
\mathcal{M}_{N,Y_i^2}^t=\int_0^t\sum_{j\in\mathscr{S}, k \in \mathscr{D}} \overline{\nabla}^{(j)}\left[\left(Y_i^{(N)}(t)\right)^2 \right]\,\widetilde{\Lambda}^\sigma_N(j,k,ds) \,, 
\end{equation}
where we have defined 
\begin{equation}\label{32CWG}
\overline{\nabla}^{(j)}\left[\left(Y_i^{(N)}(t)\right)^2 \right] :=  \left( Y_i^{(N)}(t) - j \frac{2\varphi_i(k)}{N^{\frac{3}{4}}} \right)^2 - \left(Y_i^{(N)}(t) \right)^2
\end{equation}
and
\begin{equation}\label{26CWG}
\widetilde{\Lambda}^\sigma_N(j,k,dt) := \Lambda^\sigma_N(j,k,dt) - \underbrace{N^{\frac{1}{4}} \left\vert A(j,k,N^{\frac{1}{4}}t) \right\vert  e^{-\beta j \left(N^{-\frac{1}{4}} \int \sigma d\tilde{\rho}_N(N^{\frac{1}{4}}t) + k \right)} dt}_{:= \lambda^\sigma (j,k,t)\,dt} \,.
\end{equation}
The quantity $\widetilde{\Lambda}^{\sigma}_N(j,k,dt)$ is the difference between the point process $\Lambda^{\sigma}_N(j,k,dt)$, defined on \mbox{$\mathscr{S} \times \mathscr{D} \times \mathbb{R}^+$}, and its intensity \mbox{$\lambda^{\sigma} (j,k,t)\,dt$}. 
The quantity $\left\vert A(j,k,N^{\frac{1}{4}}t) \right\vert$  indicates the number of sites $i$ that at time $N^{\frac{1}{4}} t$ have $\sigma_i = j$ and $\eta_i = k$ and it is given by
\begin{multline}\label{41CWG}
\left\vert A(j,k,N^{\frac{1}{4}}t) \right\vert = \frac{N}{4} \left[ 1 + \frac{1}{kN^{\frac{1}{4}}} \int \eta d\tilde{\rho}_N(t) + \frac{j}{N^{\frac{1}{4}}} \int \sigma d\tilde{\rho}_N(t)  \right. \\
\left. + \frac{j}{k} \left( \frac{1}{N^{\frac{1}{4}}} \int \sigma \eta d\tilde{\rho}_N(t) - \int_{\mathscr{D}} \eta \tanh(\beta \eta) \mu(d\eta) \right) \right]\,.
\end{multline}

\begin{remark}
If we call $(\mathcal{A}_t)_{t\geq0}$ the filtration generated by $\Lambda^\sigma_N$, then the processes $J_N \left[ \left(Y_i^{(N)}(t) \right)^2 \right]$ and  $\overline{\nabla}^{(j)} \left[ \left(Y_i^{(N)}(t) \right)^2 \right]$ are $\mathcal{A}_t-$adapted processes.
\end{remark}

For every index $i=1, \dots, m-1$, the following result holds. Note that it is stronger than the collapse of the processes $\left(Y_i^{(N)}\right)_{i=1}^{m-1}$, in the sense of Definition \ref{defcollapse}.

\begin{lemma}\label{33CWG}
Fix $d>2$, and assume the assumptions of Theorem \ref{CWinhom} are satisfied. Then, for every $\varepsilon >0$ there exist $N_0$ such that for every $M >0$ there is a constant $C_6>0$ for which
\begin{equation}\label{39CWG}
    \sup_{N\geq N_{0}} P\left\{\sup_{0\leq t\leq T\wedge\tau_N^M} \left( Y_i^{(N)}(t) \right)^2 > C_{6}\,N^{-\frac{1}{8} \left( 1 - \frac{2}{d} \right)}\right\}\leq\varepsilon \,.
\end{equation}
\end{lemma}

\begin{proof}
The main tool is Proposition \ref{3}. However, some assumptions in Proposition \ref{3} are not satisfied uniformly in the environment. We therefore will condition on the event
\[
A_K := \left\{ \underline{\eta} \in \mathscr{D}^N : \left| \int \sinh(\beta \eta) \varphi_i(\eta) d\hat{\rho}_N \right| + \left| \int \cosh(\beta \eta) \varphi_i(\eta) d\hat{\rho}_N \right| \leq K \right\}.
\]
The random field $\underline{\eta}$ is i.i.d., so it satisfies a standard Central Limit Theorem. Therefore, we can choose $K>0$ such that for every $N \geq 1$,
\[
P(A_K^c) \leq \frac{\varepsilon}{2}.
\]
Constants below are allowed to depend on $K$; this dependence is omitted.
We are left to show that, for every $M>0$ there is $C_6 >0$ such that 
\begin{equation}\label{39CWGcond}
    \sup_{N\geq N_{0}} P_K \left\{\sup_{0\leq t\leq T\wedge\tau_N^M} \left( Y_i^{(N)}(t) \right)^2 > C_{6}\,N^{-\frac{1}{8} \left( 1 - \frac{2}{d} \right)}\right\}\leq\frac{\varepsilon}{2} \,,
\end{equation}
where $P_K( \, \cdot \,) := P(\, \cdot \, | A_K)$. To prove (\ref{39CWGcond}) we check the conditions in Proposition \ref{3}.

\noindent
{\em Step 1}. We set $\kappa_N := N^{\frac{1}{4}}$, $\alpha_N := N^{\frac{1}{8}}$, $\beta_N \equiv 1$. Clearly
(\ref{con1}) in Proposition \ref{3} holds. 

\noindent
{\em Step 2}. We check (\ref{con2}) of Proposition \ref{3}, i.e.
\begin{equation}\label{35CWG}
    E\left[\left( Y_i^{(N)}(0) \right)^{2d}\right]\leq C_{1}\,N^{-\frac{d}{4}} \quad \mbox{for all $N$.}
\end{equation}
We start noticing that a Central Limit Theorem applies to the processes $\int \sigma \varphi_i (\eta) \, d\rho_N(0)$, since the random variables $(\sigma_j(0), \varphi_i(\eta_j))_{j=1}^N$ are independent; so, in the limit as $N \rightarrow +\infty$, $N^{\frac{1}{4}} Y_i^{(N)}(0)$ converges to a Gaussian random variable and, since $(\sigma_j(0) \varphi_i(\eta_j))_{j=1}^N$ are bounded random variables, there is convergence of all the moments. In particular (\ref{35CWG}) holds.

\noindent
{\em Step 3}. We check (\ref{con3}) of Proposition \ref{3}, i.e.
\begin{equation}\label{36CWG}
 J_N \left[ \left( Y_i^{(N)} \right)^2 \right](t) \leq -N^{\frac{1}{4}}\delta \left( Y_i^{(N)}(t) \right)^2 +C_{2} \, ,
\end{equation}
for suitable constants $\delta,C_2 >0$, which are allowed to depend on $M$, and all \mbox{$t \in [0,\tau_N^M]$} (we recall that $\beta_N \equiv 1$). Letting $X :=  - \frac{\beta}{N^{\frac{1}{4}}} \int \sigma d\tilde{\rho}_N $, we write
\[
\exp[\pm X] = 1 \pm X + R_{\pm}.
\]
Using this expansion we can perform the computation as in Proposition \ref{CWpropcritfluct}, but keeping track of the remainders  :
\begin{multline*}
\!\!\!\!\! \mbox{$\displaystyle{J_N\left[ \left( Y_i^{(N)} \right)^2 \right]}$} 
\end{multline*}
\vspace{-0.62cm}
\begin{multline*}
= N^{\frac{1}{4}}\sum_{j=1}^N \left[\cosh(\beta \eta_j) - \sigma_j \sinh(\beta \eta_j) \right] \exp\left[- \frac{\beta}{N^{\frac{1}{4}}} \! \int \sigma d\tilde{\rho}_N \right] \\
\times \left[\left( Y_i^{(N)} - \frac{2 \sigma_j}{N^{\frac{3}{4}}} \varphi_i(\eta_j)\right)^2 -  \left( Y_i^{(N)} \right)^2 \right] \end{multline*}
\vspace{-0.62cm}
\begin{multline*} 
= N^{\frac{1}{4}}\sum_{j=1}^N \left[\cosh(\beta \eta_j) - \sigma_j \sinh(\beta \eta_j) \right] \left(1-\frac{\beta \sigma_j}{N^{\frac{1}{4}}} \int \sigma d\tilde{\rho}_N + R_{\sigma_j} \right) \\ 
\times \left[ -4Y_i^{(N)} \frac{\sigma_j \varphi_i(\eta_j)}{N^{\frac{3}{4}}} + \frac{4 \varphi_i^2(\eta_j)}{N^{\frac{3}{2}}} \right] \end{multline*}
\vspace{-0.62cm}
\begin{multline}
= -4 N^{\frac{1}{4}} Y_i^{(N)} \int \sigma \mathfrak{L} \varphi_i d\tilde{\rho}_N  + 4 Y_i^{(N)}  \int \sinh(\beta \eta) \varphi_i(\eta) d\hat{\rho}_N \\ 
-  4 Y_i^{(N)}  \beta \int \sigma d\tilde{\rho}_N \int \sigma \sinh(\beta \eta) \varphi_i (\eta) d\tilde{\rho}_N \\
+  4 Y_i^{(N)} \frac{\beta}{N^{\frac{1}{4}}}\int \sigma d\tilde{\rho}_N \int \cosh(\beta \eta) \varphi_i(\eta) d\hat{\rho}_N \\
+ \frac{4}{N^{\frac{5}{4}}} \sum_{j=1}^N\left[\cosh(\beta \eta_j) + \sigma_j \sinh(\beta \eta_j)\right] \varphi_i^2(\eta_j) \\ 
+ \frac{4\beta}{N^{\frac{3}{2}}}\int \sigma d\tilde{\rho}_N \sum_{j=1}^N\left[\cosh(\beta \eta_j) + \sigma_j \sinh(\beta \eta_j)\right] \sigma_j \varphi_i^2(\eta_j) \\
+ N^{\frac{1}{4}}\sum_{j=1}^N\left[\cosh(\beta \eta_j) + \sigma_j \sinh(\beta \eta_j)\right] \left[-4 Y_i^{(N)} \frac{\sigma_j \varphi_i(\eta_j)}{N^{\frac{3}{4}}} + \frac{4 \varphi_i^2(\eta_j)}{N^{\frac{3}{2}}} \right] R_{\sigma_j}. \label{expstep3}
\end{multline}
The first term of this last expression is
\[
-4 N^{\frac{1}{4}} Y_i^{(N)} \int \sigma \mathfrak{L} \varphi_i d\tilde{\rho}_N = -4 \lambda_i N^{\frac{1}{4}} \left(Y_i^{(N)}\right)^2.
\]
We are left to show that all remaining terms are bounded, for $t \in [0,\tau_N^M]$, $\underline{\eta} \in A_K$ and assuming that in (\ref{expstep3}), $\tilde{\rho}_N$ is evaluated at time $N^{\frac{1}{4}} t$. We immediately have
\[
\left| Y_i^{(N)}(t) \right| \leq M, \ \  \left| \int \sinh(\beta \eta) \varphi_i(\eta) d\hat{\rho}_N \right| + \left| \int \cosh(\beta \eta) \varphi_i(\eta) d\hat{\rho}_N \right| \leq K.
\]
All remaining terms in (\ref{expstep3}) are of the form
\[
\int \sigma f(\eta) d\tilde{\rho}_N(N^{\frac{1}{4}} t),
\]
for some real valued $f$. Since $(\varphi_h)_{h=0}^{m-1}$ form a basis for the vector space of these functions, we can write 
\[
f = \sum_{h=0}^{m-1} \alpha_h \varphi_h.
\]
Thus
\[
\left|\int \sigma f(\eta) d\tilde{\rho}_N(N^{\frac{1}{4}} t) \right| \leq \sum_{h=0}^{m-1} \left|\alpha_h \right| \left| Y_h^{(N)}(t)\right| \leq CM,
\]
where $C$ depend on $m$, on the combinators $\alpha_h$, but not on $N$. As a consequence
\[
| R_{\pm}|  \leq \sup\left\{ e^z: |z| \leq \frac{\beta}{N^{\frac{1}{4}}} \left| \int \sigma d\tilde{\rho}_N \right| \right\} \frac{\beta^2}{2N^{\frac{1}{2}}} \left( \int \sigma d\tilde{\rho}_N \right)^2 \leq  \frac{\beta^2 M^2}{2N^{\frac{1}{2}}} e^{\beta M}.
\]
With all this, (\ref{expstep3}) implies
\[
J_N\left[ \left( Y_i^{(N)} \right)^2 \right]  \leq -4 \lambda_i N^{\frac{1}{4}} \left(Y_i^{(N)}\right)^2 + C(M)
\]
for some $M$-dependent constant $C(M)$.

\noindent
{\em Step 4}. We check (\ref{con4}) of Proposition \ref{3}, i.e. (see equation (\ref{CWmartterm}))
\begin{equation}\label{37CWG}
    \sup_{\omega\in\Omega, j\in\mathscr{S}, t\leq\tau_N^M} \left\vert \overline{\nabla}^{(j)} \left[\left( Y_i^{(N)}(t) \right)^2 \right]  \right\vert \leq C_{4}\, N^{-\frac{1}{8}} \, .
\end{equation}
For $t\leq\tau_N^M$, we easily have
\begin{align*}
  \left\vert \overline{\nabla}^{(j)}\left[ \left( Y_i^{(N)}(t) \right)^2 \right] \right\vert &=  \left\vert \left[ \frac{4\varphi_i^2(k)}{N^{\frac{3}{2}}}  - j \frac{4 \varphi_i(k) Y_i^{(N)} }{N^{\frac{3}{4}}} \right] \right\vert \\
& \leq \frac{4}{N^{\frac{3}{4}}} (1+M)  \sup_{k \in \mathscr{D}} \{ \varphi^2_i(k) + \vert \varphi_i(k) \vert \}  \leq C_4 N^{-\frac{1}{8}} \,,
\end{align*}

\noindent
{\em Step 5}. We check (\ref{con5}) of Proposition \ref{3}, i.e. (see equation (\ref{26CWG}))
\begin{equation}\label{38CWG}
    \sum_{j\in\mathscr{S},k\in\mathscr{D}}\left[ \overline{\nabla}^{(j)} \left[ \left( Y_i^{(N)}(t) \right)^2 \right] \right]^2\,\lambda^{\sigma}(j,k,t) \leq C_{5}.
\end{equation}
Recalling the definitions of $\overline{\nabla}^{(j)} \left[ \left( Y_i^{(N)} (t) \right)^2 \right]$ and $\lambda(j,k,t)$, which can be found in \eqref{32CWG} and in \eqref{26CWG}, we have
\begin{multline*}
\sum_{j\in\mathscr{S},k\in\mathscr{D}}\left[ \overline{\nabla}^{(j)} \left[ \left( Y_i^{(N)}(t) \right)^2 \right] \right]^2\,\lambda^{\sigma}(j,k,t)  \\ 
\!\!\!\!  = N^{\frac{1}{4}}  \sum_{j \in \mathscr{S}, k \in \mathscr{D}} \vert A(j,k,N^{\frac{1}{4}}t)\vert  e^{-\beta j \left(N^{-\frac{1}{4}} \int \sigma d\tilde{\rho}_N(N^{\frac{1}{4}} t) + k \right)} \\
\times \left[ \left(Y_i^{(N)}(t) -  j \frac{2\varphi_i(k)}{N^{\frac{3}{4}}} \right)^2  - \left( Y_i^{(N)}(t) \right)^2 \right]^2 
\end{multline*}
Boundedness of this last expression for $t \in [0,\tau_N^M]$, $\underline{\eta} \in A_K$ follows readily by boundedness of $Y_i^{(N)}(t)$ and $\int \sigma d\tilde{\rho}_N(N^{\frac{1}{4}} t)$ (see step 3), and the fact that \mbox{$\vert A(j,k,N^{\frac{1}{4}}t)\vert \leq N$}.

\noindent
{\em Step 6}. Conclusion. It is now enough to use (\ref{3b}). 
\end{proof}


The next step is to prove, for every $\epsilon >0$ and $N \geq 1$, the existence of a constant $M>0$ such that 
\begin{equation}\label{CWtau}
P\left\{\tau_{N}^M\leq T \right\} \leq \epsilon \,.
\end{equation}
This fact, together with Lemma \ref{33CWG}, implies the processes $Y_1^{(N)}(t),\dots,Y_{m-1}^{(N)}(t)$ converge to zero in probability, as $N$ grows to infinity, for $t$ in the whole time interval $[0,T]$. As in (\ref{39CWGcond}), we can replace $P$ by $P_K$ for a sufficiently large $K$. The idea is to consider a martingale decomposition as in (\ref{CWmartdec}) for $\psi\left(Y_0^{(N)} \right)$, where $\psi \in {\cal{C}}^1$ has bounded first derivative, and is such that $|x|>M$ implies $\psi(x)>M$; for instance, $\psi(x) = \sqrt{1+x^2}$. We obtain
\begin{equation}\label{CWmartdec1}
\psi\left(Y_0^{(N)}(t) \right) = \psi\left(Y_0^{(N)}(0) \right) + \int_0^t J_N \psi\left(Y_0^{(N)} \right) (s) ds + \mathcal{M}^t_{N,\psi},
\end{equation}
where 
\begin{equation}\label{CWmartterm1}
\mathcal{M}^t_{N,\psi} = \int_0^t \sum_{j \in \mathscr{S}, k \in \mathscr{D}}   \left[\psi \left( Y_0^{(N)}(s) - j \frac{2 \varphi_0(k)}{N^{\frac{3}{4}}} \right) - \psi\left(Y_0^{(N)}(s) \right)  \right] \widetilde{\Lambda}^\sigma_N(j,k,ds)
\end{equation}
with $\widetilde{\Lambda}^\sigma_N$ as in (\ref{26CWG}). The point now is to get bounds on $ J_N \psi\left(Y_0^{(N)}\right)$. We proceed as in (\ref{expstep3}); the only difference is in the ``gradient term'', which is now
\[
\psi \left( Y_0^{(N)}(s) - \sigma_j \frac{2 \varphi_0(\eta_j)}{N^{\frac{3}{4}}} \right) - \psi\left(Y_0^{(N)}(s) \right) = - \psi' \left(Y_0^{(N)}(s) \right)\sigma_j \frac{2 \varphi_0(\eta_j)}{N^{\frac{3}{4}}} + \mathcal{R}_N
\]
with $\mathcal{R}_N \leq \frac{C}{N^{\frac{3}{2}}}$. Proceeding as in (\ref{expstep3}), it is easily seen that
\begin{multline} \label{exptau}
\!\! J_N \psi\left(Y_0^{(N)}\right) = 2  \psi' \left(Y_0^{(N)} \right) \int \sinh(\beta \eta) \varphi_0(\eta) d\hat{\rho}_N  \\ -  2  \psi' \left(Y_0^{(N)}\right)   \beta \int \sigma d\tilde{\rho}_N \int \sigma \sinh(\beta \eta) \varphi_0 (\eta) d\tilde{\rho}_N + O_M\left(N^{\frac{1}{4}}\right),
\end{multline}
where $O_M\left(N^{\frac{1}{4}}\right)$ includes all term that, for $t \leq \tau_M^N$, are bounded by $\frac{C(M)}{N^{\frac{1}{4}}}$. The absolute value of the term
\[
2  \psi' \left(Y_0^{(N)} \right) \int \sinh(\beta \eta) \varphi_0(\eta) d\hat{\rho}_N 
\]
is bounded by $CK$, since $ \psi'$ is bounded, and $\underline{\eta} \in A_K$. For the term
\[
-  2  \psi' \left(Y_0^{(N)}\right)   \beta \int \sigma d\tilde{\rho}_N \int \sigma \sinh(\beta \eta) \varphi_0 (\eta) d\tilde{\rho}_N
\]
one should notice that $\sinh(\beta \eta) \varphi_0 (\eta)$ is orthogonal to $\varphi_0(\eta)$ in $L^2(\nu)$. This implies that $\int \sigma \sinh(\beta \eta) \varphi_0 (\eta) d\tilde{\rho}_N$ is a linear combination of $Y_1^{(N)}, Y_2^{(N)}, \ldots, Y_{m-1}^{(N)}$. Due to (\ref{39CWGcond}), we can choose a constant $C(M)$ for which this term is bounded by 
\[
C(M) N^{-\frac{1}{4} \left(1-\frac{2}{d}\right)}
\]
for $t \leq \tau_M^N$, with probability greater that $1-\frac{\varepsilon}{4}$. Denote by $B_{\varepsilon}$ the event that this bound holds true. Putting all together, we have therefore proved that in $A_K \cap B_{\varepsilon}$ and  $t \leq \tau_M^N$, 
\[
\left| J_N \psi\left(Y_0^{(N)}\right) \right| \leq CK + \frac{C(M)}{ N^{\frac{1}{4} \left(1-\frac{2}{d}\right)}}.
\]
This means that, by (\ref{CWmartdec1}), the inequality
\[
\sup_{0\leq t \leq T\wedge\tau_{N}^M} \left\vert Y_0^{(N)}(t) \right\vert \geq M
\]
implies, for $N$ and $M$ large enough, that either
\[
\left|Y_0^{(N)}(0)\right| \geq cM
\]
or
\[
\sup_{0\leq t \leq T\wedge\tau_{N}^M} \left| \mathcal{M}^t_{N,\psi} \right| \geq cM
\]
for some $c>0$. 

Thus,
\begin{allowdisplaybreaks}
\begin{align*}
\{\tau_{N}^M \leq T\}&\subseteq \bigg\{\sup_{0\leq t \leq T\wedge\tau_{N}^M} \left \{\left\vert Y_0^{(N)}(t) \right \vert,  \dots, \left\vert Y_{m-1}^{(N)}(t) \right\vert \right \} \geq M\bigg\}\\
&\subseteq \bigcup_{i=1}^{m-1} \bigg\{\sup_{0\leq t \leq T\wedge\tau_{N}^M} \left\vert Y_i^{(N)}(t) \right\vert \geq M\bigg\} \cup \left\{ \left\vert Y_0^{(N)}(0) \right \vert \geq c M \right\} \cup\\
& \hspace{100pt} \cup \bigg\{\sup_{0\leq t \leq T\wedge\tau_{N}^M} \left \vert M_{N, \psi}^t \right\vert \geq cM  \bigg\} \cup A_K^c \cup B_{\varepsilon}^c
\end{align*}
\end{allowdisplaybreaks}\\
and we obtain the following inequality for the probability of the interested set 
\begin{allowdisplaybreaks}
\begin{multline*}
P\{\tau_{N}^M \leq T\} \leq \frac{3}{4}\varepsilon + \sum_{i=1}^{m-1} P\bigg\{\sup_{0\leq t \leq T\wedge\tau_{N}^M} \left\vert Y_i^{(N)}(t) \right\vert \geq M\bigg\} \\
+ P \left\{ \left\vert Y_0^{(N)}(0) \right \vert \geq cM \right\}  + P\bigg\{\sup_{0\leq t \leq T\wedge\tau_{N}^M} \left \vert M_{N, \psi}^t \right\vert \geq cM \bigg\} \,.
\end{multline*}
\end{allowdisplaybreaks}

We estimate the three terms of the right-hand side of the inequality. 
\begin{itemize}
\item For every $i=1,\dots,m-1$, thanks to \eqref{39CWG} we have 
\[
P\bigg\{\sup_{0\leq t \leq T\wedge\tau_{N}^M} \left\vert Y_i^{(N)}(t) \right\vert \geq  M \bigg\}\leq\frac{\varepsilon}{12} \,,
\]
where for $M$ large enough.
\item Since at time \mbox{$t=0$} the spins are distributed according to a product measure, $Y_0^{(N)}(0)$ is $N^{\frac{1}{4}}$ times the sample average of independent, bounded random variables of mean zero. Therefore, for some constant $C>0$,
\[ 
E \left[ \left \vert Y_0^{(N)}(0) \right \vert \right] \leq \frac{C}{N^{\frac{1}{4}}} 
\]
and in the limit as $N \rightarrow +\infty$, we have convergence to zero in $L^1$ and then in probability. Therefore
\[ 
P \left\{ \left \vert Y_0^{(N)}(0) \right \vert \geq c M \right \} \leq \frac{\varepsilon}{12} 
\]
for $N$ sufficiently large.
\item We reduce to deal with $E \left[ \left( \mathcal{M}_{N, \psi}^T \right)^2 \right]$; in fact, by Doob's maximal inequality for martingales (we refer to Chapter VII, Section 3 of ~\cite{Shi96}) we have
\[
P \left \{\sup_{0\leq t \leq T\wedge\tau_{N}^M} \left \vert \mathcal{M}_{N, \psi}^t \right \vert \geq cM \right\} \leq \frac{E \left[ \left( \mathcal{M}_{N,\psi}^T \right)^2 \right]}{(cM)^2} \,.
\]
It is therefore enough to show that $E \left[ \left( \mathcal{M}_{N,\psi}^T \right)^2 \right]$ is bounded uniformly on $N$ and $M$. By (\ref{CWmartterm1}) and since $\psi$ is Lipschitz, we have (see also (\ref{26CWG}))
\begin{equation}\label{CWmartbound}
E \left[ \left( \mathcal{M}_{N,\psi}^T \right)^2 \right] \leq \frac{C}{N^{\frac{3}{2}}} E \left[ \int_0^T \sum_{j \in \mathscr{S}, k \in \mathscr{D}} \lambda^{\sigma} (j,k,t) dt \right].
\end{equation}
Since, by (\ref{26CWG}), $ \lambda^{\sigma} (j,k,t) \leq CN^{\frac{5}{4}}$ for some constant $C$, the boundedness of $E \left[ \left( \mathcal{M}_{N,\psi}^T \right)^2 \right]$ is established, and the proof of (\ref{CWtau}) is completed.
\end{itemize}


\subsection{Identification of the Limiting Generator and Convergence}\label{ssCW:Convergence}

We are going to show that, in the limit of infinite volume and $t \in [0,T]$, the process $Y_0^{(N)}(t)$ admits a limit in distribution, that we will be able to identify.\\

First, we need to prove the tightness of the sequence $\left\{Y_0^{(N)}(t) \right\}_{N \geq 1}$. This property implies the existence of convergent subsequences. Secondly, we will verify that all the convergent subsequences have the same limit and hence also the sequence $\left\{ Y_0^{(N)}(t) \right\}_{N \geq 1}$ must converge to that limit. 

\begin{lemma}\label{27CWG}
The sequence $\left \{ Y_0^{(N)}(t) \right\}_{N \geq 1}$ is tight.
\end{lemma}

\begin{proof}
Following \cite{CoEi88}, we use the following \emph{tightness criterion}:\\

a sequence of processes $\{\xi_N(t)\}_{N \geq 1}$ on $\mathcal{D}[0,T]$ is tight if
\begin{enumerate}
\item for every $\varepsilon>0$ there exists $M>0$ such that
\begin{equation}\label{T1}
\sup_{N}P\bigg\{\sup_{t\in[0,T]} \vert \xi_N(t)\vert \geq M\bigg\} \leq \varepsilon \,,
\end{equation}
\item for every $\varepsilon>0$ and $\alpha >0$ there exists $\delta > 0$ such that
\begin{equation}\label{T2}
\sup_{N}\sup_{0\leq \tau_1 \leq \tau_2 \leq (\tau_1 + \delta) \wedge T }P\{\vert \xi_N(\tau_2) - \xi_N(\tau_1) \vert \geq \alpha\} \leq \varepsilon \, ,
\end{equation}
where the second $\sup$ is over stopping times  $\tau_1$ and $\tau_2$, adapted to the filtration generated by the process $\xi_N$.
\end{enumerate}

We must verify the conditions \eqref{T1} and \eqref{T2} hold. Since we have already proved that for every $\epsilon>0$ the inequality $P \{ \tau_N^M \leq T \} \leq \epsilon$ is true for $M$ sufficiently large and uniformly in $N$, it is enough to show tightness for the stopped processes
\[
\left\{ Y_0^{(N)}(t \wedge \tau_N^M) \right \}_{N \geq 1} \,.
\]
We have already shown that, for $M$ large enough
\[
P\left\{ \sup_{0\leq t \leq T\wedge\tau_{N}^M} \left\vert Y_0^{(N)}(t) \right\vert \geq M \right\} \leq \varepsilon
\]
which yields \eqref{T1}. To obtain \eqref{T2}, we notice that
\begin{equation}\label{CWti1}
\left\vert Y_0^{(N)}(\tau_2)  -  Y_0^{(N)}(\tau_1) \right\vert = \left\vert \int_{\tau_1}^{\tau_2} J_N \left( Y_0^{(N)}(u) \right) du + \mathcal{M}_{N, Y_0 }^{\tau_1,\tau_2} \right\vert \,,
\end{equation}
where we have denoted
\[
\mathcal{M}_{N, Y_0 }^{\tau_1,\tau_2} = -\frac{2}{N^{\frac{3}{4}}}\int_{\tau_1}^{\tau_2} \sum_{j \in \mathscr{S}, k \in \mathscr{D}} j \varphi_0(k) \, \widetilde{\Lambda}_N^\sigma(j,k,du) 
\]
and $\widetilde{\Lambda}_N^\sigma$ is as in \eqref{26CWG}. As in the proof of Lemma \ref{33CWG}, one shows that both $ J_N \left( Y_0^{(N)}\right)$ and the quadratic variation of $\mathcal{M}_{N, Y_0 }^{\tau_1,\tau_2} $ are uniformly bounded in $N$, from which \eqref{T2} follows for the processes $\left\{ Y_0^{(N)}(t \wedge \tau_N^M) \right \}_{N \geq 1}$.
\end{proof}

Lemma \ref{27CWG} implies that there exist convergent subsequences for the sequence $\left\{ Y_0^{(N)}(t) \right\}_{N \geq 1}$. With abuse of notation, let $\left\{ Y_0^{(n)}(t)  \right\}_{n \geq 1}$ denote one of such a subsequence and let $\psi \in \mathcal{C}_b^2$. The following decomposition holds
\begin{equation}\label{28CWG}
\psi \left( Y_0^{(n)}(t) \right) - \psi \left( Y_0^{(n)}(0) \right) = \int_0^t J_n \psi \left( Y_0^{(n)}(u) \right)du + \mathcal{M}_{n,\psi}^t \,,
\end{equation}
where
\begin{multline*}
J_n \psi \left( Y_0^{(n)}(t) \right) =2 \psi'\left( Y_0^{(n)}(t) \right) \left\{ n^{\frac{1}{4}} \int \tanh(\beta\eta) d\tilde{\rho}_n(t)  \right. \\
+ \left. \beta \int \sigma d\tilde{\rho}_n(t) \int d\tilde{\rho}_n(t) - \beta \int \sigma d\tilde{\rho}_n(t) \int \sigma \tanh(\beta \eta) d\tilde{\rho}_n(t) \right\} + o_M(1).
\end{multline*}
The remainder $o_M(1)$ goes to zero as $n \rightarrow +\infty$, for $t \leq \tau_n^M$. 
If we compute the limit as $n\rightarrow +\infty$, using the facts that a Central Limit Theorem applies to the term $\int \tanh(\beta\eta) d\tilde{\rho}_n(t)$, the integral $\int d\tilde{\rho}_n(t)$ is zero since $\widetilde{\rho}_n$ is a centered measure, and the process $\int \sigma \tanh(\beta \eta) d\tilde{\rho}_n(t)$ collapse since $\tanh(\beta \eta)$ and $\varphi_0(\eta) = \frac{1}{\cosh(\beta \eta)}$ are orthogonal in $L^2 \left( \nu \right)$,  we have, in the sense of weak convergence of processes: 
\[ 
J_n \psi\left( Y_0^{(n)}(t \wedge \tau_n^M) \right) \xrightarrow[\quad w \quad]{n\rightarrow +\infty} J \psi(Y_0(t\wedge \tau_n^M)) 
\quad \mbox{with} \quad
J \psi(Y_0(t)) = 2 \, \mathscr{H} \,  \psi'(Y_0(t)) 
\]
and where $\mathscr{H}$ is a Gaussian random variable. Then, because of \eqref{28CWG} and \eqref{CWtau}, we obtain 
\[  
\mathcal{M}_{n,\psi}^t \xrightarrow[\quad w \quad]{n \rightarrow +\infty}  \mathcal{M}_{\psi}^t := \psi(Y_0(t)) - \psi(Y_0(0)) - \int_0^t J \psi(Y_0(u))du \, ,\] 
for $t \in [0,T]$.
We must prove the following Lemma:

\begin{lemma}\label{31CWG}
$\mathcal{M}_\psi^t$ is a martingale (with respect to $t$); in other words, for all $s,t \in [0,T]$, $s\leq t$ and for all measurable and bounded functions  $g(Y_0([0,s]))$ the following identity holds: 
\begin{equation}\label{29CWG}
E[ \mathcal{M}_{\psi}^t g(Y_0([0,s]))] = E[ \mathcal{M}_{\psi}^s g(Y_0([0,s]))] \,.
\end{equation}
\end{lemma}

\begin{proof}
We begin by showing that \eqref{29CWG} follows from the fact, that will be proved later, that for every $t$ fixed,
$ \{ \mathcal{M}_{n,\psi}^t \}_{n \geq 1}$ is a uniformly integrable sequence of random variables.\\
Since $\mathcal{M}_{n,\psi}^t$ is a martingale (with respect to $t$) for every $n$, we have that for all $s,t \in [0,T]$, $s\leq t$ and for all measurable and bounded functions $g(Y_0([0,s]))$ 
\begin{equation}\label{CWnmartlim}
E[ \mathcal{M}_{n,\psi}^t g(Y_0([0,s]))] = E[ \mathcal{M}_{n,\psi}^s g(Y_0([0,s]))]. 
\end{equation}
Now, as we have seen, $\mathcal{M}_{n,\psi}^t$ and $\mathcal{M}_{n,\psi}^s$ have a weak limit; this, together with uniform integrability, imply convergence in $L^1$. Thus \eqref{29CWG} follows by taking limit in \eqref{CWnmartlim}.
%
%
%
\\
It remains to check  that $\{ \mathcal{M}_{n,\psi}^t \}_{n \geq 1}$ is a uniformly integrable family. A sufficient condition for uniform integrability is that $\sup_n E[\vert  \mathcal{M}_{n,\psi}^t \vert ^2]< +\infty$ (see again ~\cite{Shi96}). \\ This, however, is exactly what we have done already in \eqref{CWmartbound}.
\end{proof}

\paragraph{Proof of Theorem \ref{CWinhom}.}
We have shown that any weak limit of $Y_0^{(n)}(\, \cdot \,)$ solves the martingale problem with infinitesimal generator $J$, which admits a unique solution. It follows that all convergent subsequences have the same limit and so the sequence itself converges to that limit.

\section{Proofs for the Random Kuramoto Model}\label{proofsK}

Throughout this section we assume $\omega \leq \frac{1}{2 \sqrt{2}}$, even though this assumption will be relevant only starting from Section \ref{ssK:Collapse}. Whenever needed, we will comment on the necessary changes to cover the case $ \frac{1}{2 \sqrt{2}} < \omega < \frac{1}{2}$.

\subsection{Preliminaries}

\paragraph{Proof of Lemma \ref{lmm:specK}.} If $\varphi(\cdot,\cdot)$ belongs to the null space of $\mathfrak{L}$, then $\mathfrak{L} \varphi=0$. Therefore, we require that
\begin{multline}\label{31Ki}
\frac{1}{2} \frac{\partial^2 \varphi}{\partial x^2}(x,\eta) + \omega\eta \frac{\partial \varphi}{\partial x}(x,\eta) + (1+4\omega^2) \bigg[ \cos x \, \frac{1}{2\pi} \int \cos y \, \varphi(y,\varsigma) \, q_*(dy,d\varsigma)  \\
+ \sin x \, \frac{1}{2\pi} \int \sin y \, \varphi(y,\varsigma) \, q_*(dy,d\varsigma) \bigg] =0 \,.
\end{multline}
We solve the ordinary differential equation \eqref{31Ki}. Having defined
\begin{equation}\label{32Ki}
A:=\frac{1}{2\pi} \int \cos y \, \varphi(y,\varsigma) \, q_*(dy,d\varsigma) \mbox{ and } B:=\frac{1}{2\pi} \int \sin y \, \varphi(y,\varsigma) \, q_*(dy,d\varsigma) \,,
\end{equation}
the solution is $\varphi(x,\eta) = 2(B - 2 A \omega \eta) \sin x + 2(A + 2 B \omega \eta ) \cos x$; this function yields a solution of \eqref{31Ki} provided that it satisfies the self-consistency relations \eqref{32Ki}, but it does for every value of $A$ and $B$. Then the two directions which generate the kernel are  $\sin x + 2 \omega \eta \cos x$ and $\cos x - 2\omega\eta \sin x$.

\begin{remark}
In the case that $\theta \neq 1 + 4\omega^2$, the unique value for which the self-consistency relations in \eqref{32Ki} are satisfied is $A=B=0$, meaning that at the critical point the kernel of the operator $\mathfrak{L}$ is two-dimensional, while it is trivial  for all the other values of the parameter $\theta$. 
\end{remark}

The part of the statement of Lemma \ref{lmm:specK} concerning spectrum and eigenspaces is easily proved by direct computation, and the fact that the set $\{v_k^{(i)}: k \geq 1, i=1,2,3,4\}$ spans a dense subset of $L^2([0,2\pi) \times \{-1,1\})$.
\subsection{Perturbation Theory}\label{ssK:PerturbationTheory}

In the rest of the section, we often consider the \emph{time-rescaled} infinitesimal generator $J_N = \sqrt{N} L_N$, where $L_N$ is given by \eqref{Kfluctuations}. To determine the limiting generator $J$, we need to apply the first order perturbation theory. The methodology for treating a perturbation problem has been developed in the paper \cite{PaStVa77} and extends the earlier works done in \cite{Krt73, Pap77}.

It will be useful to keep in mind the following simple fact, which is just a restatement of Proposition \ref{Kpropfluct}.

\begin{proposition}\label{Kpropcritfluct}
Under the assumptions of Proposition \ref{Kpropfluct}, we have
\begin{equation} \label{Krescaledfluctuations}
J_N \psi \left( \int \phi(x,\eta) d\tilde{\rho}_N \right) =  \sqrt{N} L^{(1)} \psi + N^{\frac{1}{4}} L^{(2)} \psi +  L^{(3)} \psi + N^{-\frac{1}{4}} L^{(4)} \psi
\end{equation}
where
\begin{align*}
L^{(1)} \psi &:= \sum_{i=1}^n \partial_i \psi \left( \int \phi(x, \eta) d\tilde{\rho}_N \right) \int \mathfrak{L} \phi_i(x, \eta) d\tilde{\rho}_N \\
L^{(2)} \psi &:= \theta \sum_{i=1}^n \partial_i \psi \left( \int \phi(x, \eta) d\tilde{\rho}_N \right) \int \frac{\partial \phi_i}{\partial x} (x, \eta) \sin(y-x) d\tilde{\rho}_N d\tilde{\rho}_N \\
L^{(3)} \psi &:= \frac{1}{2} \sum_{i,k=1}^{n} \partial_{ik}^2 \psi \left( \int \phi(x, \eta) d\tilde{\rho}_N \right)  \int \frac{\partial \phi_i}{\partial x} (x,\eta) \frac{\partial \phi_k}{\partial x} (x,\eta) dq_* \\
L^{(4)} \psi &:= \frac{1}{2} \sum_{i,k=1}^{n} \partial_{ik}^2 \psi \left( \int \phi(x, \eta) d\tilde{\rho}_N \right) \int \frac{\partial \phi_i}{\partial x} (x,\eta) \frac{\partial \phi_k}{\partial x} (x,\eta) d\tilde{\rho}_N \,.
\end{align*}
\end{proposition}

As first step (Section \ref{ssK:Collapse}) we show that for every $\phi \in {\mbox span}\left\{v_1^{(3)}, v_1^{(4)}, v_k^{(i)}: k \geq 2, i=1,2,3,4 \right\}$, the process
\[
\int \phi (x, \eta) \, d\tilde{\rho}_N(\sqrt{N}t)
\]
collapses to zero in the sense of Definition \ref{defcollapse}. We are therefore left to understand the behavior as $N \rightarrow +\infty$ of the two-dimensional process
\[
\left(V_1^{(1,N)}(t), V_1^{(2,N)}(t) \right) := \left(\int v_1^{(1)} d\tilde{\rho}_N(\sqrt{N}t), \int v_1^{(2)} d\tilde{\rho}_N(\sqrt{N}t) \right).
\]
For this reason, for $\psi \in \mathcal{C}^2(\mathbb{R}^2,\mathbb{R})$, we need to control 
\[
J_N \psi\left(\int v_1^{(1)} d\tilde{\rho}_N, \int v_1^{(2)} d\tilde{\rho}_N \right).
\]
The first term in the r.h.s. of \eqref{Krescaledfluctuations} vanishes, since $v_1^{(1)}, v_1^{(2)} \in \ker(\mathfrak{L})$. In order to compensate for the second diverging term $ N^{\frac{1}{4}} L^{(2)} \psi$, one introduces a ``small'' perturbation of $\psi$ of the form
\begin{equation} \label{Kperturb}
\psi_N = \psi + N^{-\frac{1}{4}} \psi_1,
\end{equation}
for some $\psi_1$ to be chosen. We obtain
\begin{equation} \label{Kexpansion}
J_N \psi_N = N^{\frac{1}{4}}\left[ L^{(2)} \psi + L^{(1)} \psi_1 \right] + L^{(3)} \psi + L^{(2)} \psi_1 + o(1).
\end{equation}
In order to avoid divergence, $\psi_1$ should be chosen in such a way that $L^{(2)} \psi + L^{(1)} \psi_1 = 0$. At a purely formal level we are led to set
\begin{equation} \label{psi1}
 \psi_1 := - \left(L^{(1)}\right)^{-1} L^{(2)} \psi \,,
 \end{equation}
which gives
\begin{equation} \label{Klimgen}
 J_N \psi_N \, \xrightarrow{N \rightarrow +\infty} \, \left[L^{(3)} - L^{(2)} \left(L^{(1)}\right)^{-1} L^{(2)} \right] \psi =: J \psi.
 \end{equation}
The operator $J$ is therefore the candidate for the generator of the limiting process $\left( V^{(1)}, V^{(2)} \right)$. In order to make a rigorous proof out of this formal argument, the following two steps are needed:
\begin{enumerate}
\item
The operator $ \left(L^{(1)}\right)^{-1} L^{(2)}$ has to be properly defined.
\item
From the above convergence of operators one must derive weak convergence of processes. 
\end{enumerate}
Step 2 will be dealt with in Section \ref{ssK:Convergence}, through standard martingale techniques. We consider now step 1. The needed computations are rather long, but follow few basic ideas, that we now illustrate. First observe that
\begin{equation} \label{pertth1}
 L^{(2)} \psi\left(\int v_1^{(1)} d\tilde{\rho}_N, \int v_1^{(2)} d\tilde{\rho}_N \right) = \theta \sum_{i=1}^2 \partial_i \psi(\cdot \, , \, \cdot) \int \frac{\partial v_1^{(i)}}{\partial x} (x,\eta) \sin(y-x) d\tilde{\rho}_N d\tilde{\rho}_N.
 \end{equation}
We give the details for the term  $\int \frac{\partial v_1^{(1)}}{\partial x} (x,\eta) \sin(y-x) d\tilde{\rho}_N d\tilde{\rho}_N$, the other being similar. Letting
\[
V_k^{(i,N)} := \int v_k^{(i)}  d\tilde{\rho}_N,
\]
by applying standard trigonometric formulas we obtain
\begin{multline} \label{pertth2}
\int \frac{\partial v_1^{(1)}}{\partial x} (x,\eta) \sin(y-x) d\tilde{\rho}_N d\tilde{\rho}_N \\ = \frac{1}{4(1-4\omega^2)} \left[\left(V_1^{(1,N)} - 2 \omega V_1^{(4,N)}\right) \left(V_2^{(2,N)}+ V_2^{(4,N)}\right) \right. \\ + \left. \left(V_1^{(2,N)} - 2 \omega V_1^{(3,N)}\right) \left(V_2^{(1,N)}+ V_2^{(3,N)}\right)\right] \\ - \frac{i \omega}{2 (1-4\omega^2)} \left[\left(V_1^{(1,N)} - 2 \omega V_1^{(4,N)}\right) \left(V_2^{(4,N)}- V_2^{(2,N)}\right) \right. \\ + \left. \left(V_1^{(2,N)} - 2 \omega V_1^{(3,N)}\right) \left(V_2^{(1,N)}- V_2^{(3,N)}\right)\right] .
\end{multline}
This means that $L^{(2)} \psi\left(\int v_1^{(1)} d\tilde{\rho}_N, \int v_1^{(2)} d\tilde{\rho}_N \right)$ is a linear combination of terms of the form
\begin{equation} \label{pertth3}
\partial_i \psi \left(\int v_1^{(1)} d\tilde{\rho}_N, \int v_1^{(2)} d\tilde{\rho}_N \right) \int v_1^{(j)} d\tilde{\rho}_N  \int v_2^{(h)} d\tilde{\rho}_N ,
\end{equation}
$i=1,2$, $j,h = 1,2,3,4$. If we denote by $\lambda_k^j$ the eigenvalue of $\mathfrak{L}$ corresponding to the eigenfunction $v_k^{(j)}$, in the critical case $\theta = 1 + \omega^2$, we easily obtain
\begin{multline*}
L^{(1)} \left[\frac{1}{\lambda_1^j + \lambda_2^h} \partial_i \psi \left(\int v_1^{(1)} d\tilde{\rho}_N, \int v_1^{(2)} d\tilde{\rho}_N \right) \int v_1^{(j)} d\tilde{\rho}_N  \int v_2^{(h)} d\tilde{\rho}_N \right]  \\ =  \partial_i \psi \left(\int v_1^{(1)} d\tilde{\rho}_N, \int v_1^{(2)} d\tilde{\rho}_N \right) \int v_1^{(j)} d\tilde{\rho}_N  \int v_2^{(h)} d\tilde{\rho}_N;
\end{multline*}
this defines $\left(L^{(1)}\right)^{-1}$ for the whole expression in \eqref{pertth2}. Thus, the perturbation \eqref{psi1} is now well defined.  

\noindent
A further comment is relevant. In the expression for the limiting generator in \eqref{Klimgen}, the quantity
\[
 L^{(2)}\left(L^{(1)}\right)^{-1} L^{(2)}\psi 
 \]
 appears. Moreover, we have seen that $\left(L^{(1)}\right)^{-1} L^{(2)}\psi $ is linear combination of terms as in \eqref{pertth3}. We will prove later that, when evaluated at time $\sqrt{N} t$,
 \begin{itemize}
 \item
 the sequences of  processes $ \int v_1^{(j)} d\tilde{\rho}_N$, $j=3,4$,  and $ \int v_2^{(h)} d\tilde{\rho}_N$, $h=1,2,3,4$ collapse to zero;
 \item
  the sequences of processes $ \int v_1^{(j)} d\tilde{\rho}_N$, $j=1,2$ are {\em tight}.
  \end{itemize}
In particular, the processes  $\left(L^{(1)}\right)^{-1} L^{(2)}\psi$  collapse to zero. We then have to apply $L^{(2)}$ again. It is easy to show what follows.
\begin{itemize}
 \item
When $\partial_i \psi \left(\int v_1^{(1)} d\tilde{\rho}_N, \int v_1^{(2)} d\tilde{\rho}_N \right) \int v_1^{(j)} d\tilde{\rho}_N  \int v_2^{(h)} d\tilde{\rho}_N$ has $j=3,4$, i.e. it has ``two collapsing factors'', then 
\[
L^{(2)} \left[\partial_i \psi \left(\int v_1^{(1)} d\tilde{\rho}_N, \int v_1^{(2)} d\tilde{\rho}_N \right) \int v_1^{(j)} d\tilde{\rho}_N  \int v_2^{(h)} d\tilde{\rho}_N \right]
\]
is still collapsing to zero.
\item
When $j=1,2$, non collapsing terms in the expression above, arise from
\[
\partial_i \psi \left(\int v_1^{(1)} d\tilde{\rho}_N, \int v_1^{(2)} d\tilde{\rho}_N \right) \int v_1^{(j)} d\tilde{\rho}_N \int \frac{\partial v_2^{(h)}}{\partial x} (x,\eta) \sin(y-x)  d\tilde{\rho}_N  d\tilde{\rho}_N,
\]
since when the Prostapheresis formulas are applied to $\frac{\partial v_2^{(h)}}{\partial x} (x,\eta) \sin(y-x)$, terms of the form $ \int v_1^{(j)} d\tilde{\rho}_N  \int v_1^{(l)} d\tilde{\rho}_N$, $j,l=1,2$, appear. 
\end{itemize}
Carefully performing a long but straightforward calculation, one obtains the following statement.
\begin{proposition} \label{Kpropperthth}
Up to collapsing terms (as the symbol $\simeq$ is intended to mean) we have
\begin{multline*}
 L^{(2)}\left(L^{(1)}\right)^{-1} L^{(2)}\psi \left(\int v_1^{(1)} d\tilde{\rho}_N, \int v_1^{(2)} d\tilde{\rho}_N \right)  \\ \simeq - \frac{(1+4\omega^2)^2(1-8\omega^2)}{4(1-4\omega^2)^3(1+\omega^2)} V_1^{(1,N)}(t) \left[ \left( V_1^{(1,N)}(t) \right)^2 + \left( V_1^{(2,N)}(t) \right)^2 \right] \partial_1 \psi(\cdot, \cdot) \\ - \frac{(1+4\omega^2)^2(1-8\omega^2)}{4(1-4\omega^2)^3(1+\omega^2)} V_1^{(2,N)}(t) \left[ \left( V_1^{(1,N)}(t) \right)^2 + \left( V_1^{(2,N)}(t) \right)^2 \right] \partial_2 \psi(\cdot, \cdot),
\end{multline*}
with
\[
V_k^{(i,N)} := \int v_k^{(i)}  d\tilde{\rho}_N \,.
\]
 
 \end{proposition}

\subsection{Collapsing Processes}\label{ssK:Collapse}

From now on we always assume $\theta = 1+4\omega^2$, with $\omega < \frac{1}{2}$. 
In what follows, it is more convenient to work with the following real-valued basis of $L^2([0,2\pi) \times \{-1,1\})$:
\[
\left\{ v_1^{(i)}, y_h^{(i)}: i=1,2,3,4, \, h \geq 2 \right\},
\]
where 
\[
y_h^{(1)}(x,\eta) := \cos hx \ \ y_h^{(2)}(x,\eta) := \sin hx \ \ y_h^{(3)}(x,\eta) := \eta \cos hx  \ \ y_h^{(4)}(x,\eta) := \eta \sin hx.
\]
We also set
\[
Y_h^{(i,N)} := \int y_h^{(i)} d\tilde{\rho}_N,
\]
and write $Y_h^{(i,N)}(t)$ for $\int y_h^{(i)} d\tilde{\rho}_N(\sqrt{N} t)$. \\
For $r \geq 1$ define
\begin{equation} \label{Knormr}
\left\| \tilde{\rho}_N \right\|^2_r := \left(V_1^{(3,N)}\right)^2 + \left(V_1^{(4,N)}\right)^2 + \sum_{i=1}^4 \sum_{h \geq 2} \frac{1}{\left(1+h^2\right)^r} \left(Y_h^{(i,N)}\right)^2.
\end{equation}
Clearly, showing that  the sequences of  processes $ \int v_1^{(j)} d\tilde{\rho}_N$, $j=3,4$,  and $ \int v_2^{(h)} d\tilde{\rho}_N$, $h=1,2,3,4$ collapse to zero is equivalent to show that   the sequences of  processes $ \int v_1^{(j)} d\tilde{\rho}_N$, $j=3,4$,  and $ \int y_2^{(h)} d\tilde{\rho}_N$, $h=1,2,3,4$ collapse to zero which, in turn, is implied by the fact that the sequence $\left\| \tilde{\rho}_N \right\|^2_r $ collapses to zero. All processes here are meant to be evaluated at time $\sqrt{N} t$.
For $N \geq 1$, $M >0$ define
\[
\tau_N^M := \inf_{t \geq 0} \left\{ \left\| \tilde{\rho}_N(\sqrt{N} t) \right\|^2_r \geq M  \mbox{ or } \left|V_1^{(1,N)}(t)\right| \geq M \mbox{ or } \left|V_1^{(2,N)}(t)\right| \geq M\right\}.
\]

Our first result concerns collapsing of the stopped process $\left\| \tilde{\rho}_N(\sqrt{N}(t \wedge \tau_N^M) \right\|^2_r $.
\begin{lemma}\label{12Kinh}
Fix $d>2$ and $r>\frac{3}{2}$. Then, for every $\varepsilon>0$ and $M>0$, there exist $N_0 > 0$ and $C_5>0$,  for which
\begin{equation}\label{18Kinh}
 \sup_{N\geq N_{0}} P\left\{\sup_{0\leq t\leq T\wedge\tau_N^M}\| \tilde{\rho}_N(\sqrt{N}t) \|^2_r > C_{5} \, N^{\frac{1}{2d} - \frac{1}{4}} \right\}\leq\varepsilon \,.
\end{equation}
\end{lemma}

\begin{proof}
We apply Proposition \ref{3cont}. We set $\kappa_N = \sqrt{N}$, $\alpha_N = N^{\frac{1}{4}}$, $\beta_N = N^{\frac{1}{4}}$. Conditions \eqref{con1cont} and \eqref{con2cont} of Proposition \ref{3cont} are easy to check. We are therefore left to check conditions \eqref{con3cont} and \eqref{con5cont}. We observe that $\left\| \tilde{\rho}_N(\sqrt{N} t) \right\|^2_r$ admits the semimartingale representation
\[
d\left\| \tilde{\rho}_N(\sqrt{N} t) \right\|^2_r = J_N \left\| \tilde{\rho}_N \right\|^2_r (\sqrt{N} t) dt + N^{\frac{1}{4}}\sum_{j=1}^N \frac{\partial}{\partial x_j} \left\| \tilde{\rho}_N \right\|^2_r dW_j(t),
\]
where $\{W_j(t) : t>0, j=1,\dots, N\}$ is a system of independent standard Brownian motions on $[0,2\pi]$. We show the following inequalities for every $t \in [0,\tau_N^M]$, which imply \eqref{con3cont} and \eqref{con5cont}:
\begin{equation} \label{Kcollaps1}
J_N \left\| \tilde{\rho}_N \right\|^2_r (\sqrt{N} t) \leq - \left(\frac{1}{2} - 2 \omega^2\right) \sqrt{N} \left\| \tilde{\rho}_N \right\|^2_r (\sqrt{N} t) + C N^{\frac{1}{4}},
\end{equation}
\begin{equation} \label{Kcollaps2}
\sqrt{N} \sum_{j=1}^N \left(\frac{\partial}{\partial x_j} \left\| \tilde{\rho}_N\right\|^2_r\right) ^2 \leq C
\end{equation}
for some constant $C$, that is allowed to depend on $M$.

\noindent
{\em Step 1: proof of \eqref{Kcollaps1}}. We use \eqref{Krescaledfluctuations}:
\begin{equation} \label{Kcollaps3}
J_N \left\| \tilde{\rho}_N \right\|^2_r =  \sqrt{N} L^{(1)} \left\| \tilde{\rho}_N \right\|^2_r  + N^{\frac{1}{4}} L^{(2)} \left\| \tilde{\rho}_N \right\|^2_r  +  L^{(3)} \left\| \tilde{\rho}_N \right\|^2_r  + N^{-\frac{1}{4}} L^{(4)} \left\| \tilde{\rho}_N \right\|^2_r.
\end{equation}
We begin to deal with $L^{(1)} \left\| \tilde{\rho}_N \right\|^2_r$. Due to uniform convergence of the series defining  $\left\| \tilde{\rho}_N \right\|^2_r$, we can apply $L^{(1)} $ term by term. For $i=3,4$
\begin{equation} \label{Kcollaps4}
L^{(1)}\left( \int v_1^{(i)} d\tilde{\rho}_N \right)^2 = - \left(1- 4 \omega^2\right)\left( \int v_1^{(i)} d\tilde{\rho}_N \right)^2.
\end{equation}
Also, by direct computation,
\begin{equation} \label{Kcollaps5}
L^{(1)} \sum_{i=1}^4 \left( \int y_h^{(i)} d\tilde{\rho}_N \right)^2 = -2h^2  \sum_{i=1}^4 \left( \int y_h^{(i)} d\tilde{\rho}_N \right)^2.
\end{equation}
Letting $\lambda := 1 - 4 \omega^2 >0$, by \eqref{Kcollaps4} and \eqref{Kcollaps5} we obtain
\begin{equation} \label{Kcollaps6}
L^{(1)} \left\| \tilde{\rho}_N \right\|^2_r 
= - \lambda \left\| \tilde{\rho}_N \right\|^2_r 
+ \sum_{h \geq 2} \frac{\lambda - 2 h^2}{\left(1+h^2\right)^r}
\sum_{i=1}^4 \left( \int y_h^{(i)} d\tilde{\rho}_N \right)^2.
\end{equation}
We now compute $L^{(2)} \left\| \tilde{\rho}_N \right\|^2_r$. A ``typical'' summand with $h \geq 2$ of the infinite sum giving $L^{(2)} \left\| \tilde{\rho}_N \right\|^2_r$,  is
\[
L^{(2)} \left( \int y_h^{(i)} d\tilde{\rho}_N \right)^2 = 2(1+4\omega^2) \int y_h^{(i)} d\tilde{\rho}_N \int \frac{\partial y_h^{(i)}}{\partial x} \sin(y-x) d\tilde{\rho}_N d\tilde{\rho}_N.
\]
By using Prostapheresis formulas, one realizes that $\int y_h^{(i)} d\tilde{\rho}_N \int \frac{\partial y_h^{(i)}}{\partial x} \sin(y-x) d\tilde{\rho}_N d\tilde{\rho}_N$ is a linear combination, with uniformly bounded coefficients, of terms of the form
\begin{equation}\label{Kkeypoint}
Y_1^{(j,N)} Y_h^{(i,N)} Y_{h \pm 1}^{(l,N)} \leq Y_1^{(j,N)} \left[\left(Y_h^{(i,N)}\right)^2 + \left(Y_{h \pm 1}^{(l,N)}\right)^2\right].
\end{equation}
Summing over $h \geq 2$ and observing that, for $t \in [0,\tau_N^M]$, $Y_1^{(j,N)}(t) \leq cM$ for some constant $c$, we obtain (omitting the evaluation at $\sqrt{N} t$)
\begin{equation} \label{Kcollaps7}
L^{(2)}\left[\sum_{h \geq 2} \frac{1}{\left(1+h^2\right)^r}\sum_{i=1}^4 \left( \int y_h^{(i)} d\tilde{\rho}_N \right)^2\right] \leq C(M) \sum_{h \geq 2} \frac{1}{\left(1+h^2\right)^r}\sum_{i=1}^4
\left( \int y_h^{(i)} d\tilde{\rho}_N \right)^2,
\end{equation}
for some $M$-dependent constant $C(M)$. As far as the first two summands of $\left\| \tilde{\rho}_N \right\|^2_r$ are concerned by similar arguments a rough bound for $t \in [0,\tau_N^M]$ of the form
\begin{equation} \label{Kcollaps8}
L^{(2)}\left[\sum_{i=3,4} \left(\int v_1^{(i)} d\tilde{\rho}_N\right)^2 \right] \leq C(M)
\end{equation}
is obtained. Putting together \eqref{Kcollaps6}, \eqref{Kcollaps7} and \eqref{Kcollaps8}, 
\begin{align} \label{Kcollaps9}
\sqrt{N} & L^{(1)} \left\| \tilde{\rho}_N \right\|^2_r  + N^{\frac{1}{4}} L^{(2)} \left\| \tilde{\rho}_N \right\|^2_r \nonumber\\ 
& \leq  - \sqrt{N} \lambda \left\| \tilde{\rho}_N \right\|^2_r + \sqrt{N} \sum_{h \geq 2} \frac{\lambda - 2 h^2}{\left(1+h^2\right)^r} \sum_{i=1}^4 \left( \int y_h^{(i)} d\tilde{\rho}_N \right)^2 \nonumber\\ 
& \quad + N^{\frac{1}{4}}C(M) \sum_{h \geq 2} \frac{1}{\left(1+h^2\right)^r}\sum_{i=1}^4 \left( \int y_h^{(i)} d\tilde{\rho}_N \right)^2 + N^{\frac{1}{4}} C(M) \nonumber\\ 
& \leq  - \sqrt{N} \lambda \left\| \tilde{\rho}_N \right\|^2_r + N^{\frac{1}{4}} C(M)
\end{align}
where the last inequality holds for $N$ sufficiently large so that
\[
\sqrt{N} \left(\lambda - 2 h^2\right) + N^{\frac{1}{4}}C(M) \leq 0
\]
for every $h \geq 2$. \\ Consider now the term $L^{(3)} \left\| \tilde{\rho}_N \right\|^2_r$. We have
\begin{multline} \label{Kcollaps10}
L^{(3)} \left\| \tilde{\rho}_N \right\|^2_r = 2 \int \left( \frac{\partial v_1^{(3)}}{\partial x} \right)^2 dq_* + 2 \int \left( \frac{\partial v_1^{(4)}}{\partial x} \right)^2 dq_* \\ + 2 \sum_{h \geq 2}\frac{1}{\left(1+h^2\right)^r} \sum_{i=1}^4 \int \left( \frac{\partial y_h^{(i)}}{\partial x} \right)^2 dq_*.
\end{multline}
By the simple bound
\[
\left( \frac{\partial y_h^{(i)}}{\partial x} \right)^2 \leq h^2,
\]
using the fact that for $r>\frac{3}{2}$ 
\[
\sum_{h \geq 2} \frac{h^2}{\left(1+h^2\right)^r} < +\infty,
\]
from \eqref{Kcollaps10} we get
\begin{equation} \label{Kcollaps11}
L^{(3)} \left\| \tilde{\rho}_N \right\|^2_r \leq C
\end{equation}
for some constant $C$. The treatment of the term $L^{(4)} \left\| \tilde{\rho}_N \right\|^2_r$ is quite similar, since it is obtained from \eqref{Kcollaps10} replacing $q_*$ with $\tilde{\rho}_N$. Having $\tilde{\rho}_N$ total variation $N^{\frac{1}{4}}$, we get
\begin{equation} \label{Kcollaps12}
L^{(4)} \left\| \tilde{\rho}_N \right\|^2_r \leq CN^{\frac{1}{4}}.
\end{equation}
By \eqref{Kcollaps9}, \eqref{Kcollaps11} and \eqref{Kcollaps12}, \eqref{Kcollaps1} follows.

\noindent
{\em Step 2: proof of \eqref{Kcollaps2}}. Consider the summand 
\[
\left( \int y_{h}^{(i)} d\tilde{\rho}_N \right)^2
\]
of $\left\| \tilde{\rho}_N \right\|^2_r$. The summands containing $v_1^{(i)}$ are dealt with similarly. We have
\[
\frac{\partial}{\partial x_j} \left( \int y_{h}^{(i)} d\tilde{\rho}_N \right)^2 = \frac{2}{N^{\frac{3}{4}}} \left( \int y_{h}^{(i)} d\tilde{\rho}_N \right) \frac{\partial  y_{h}^{(i)}}{\partial x} (x_j, \eta_j),
\]
so that 
\[
\left| \frac{\partial}{\partial x_j} \left( \int y_{h}^{(i)} d\tilde{\rho}_N \right)^2 \right| \leq \frac{2 h}{N^{\frac{3}{4}}} \left| \int y_{h}^{(i)} d\tilde{\rho}_N \right|.
\]
Thus
\begin{multline*}
\left| \frac{\partial}{\partial x_j} \sum_{h \geq 2} \frac{1}{\left(1+h^2\right)^r} \sum_{i=1}^4  \left( \int y_{h}^{(i)} d\tilde{\rho}_N \right)^2 \right|  
\\ \leq \frac{2}{N^{\frac{3}{4}}} \sum_{h \geq 2} \frac{h}{\left(1+h^2\right)^r} \sum_{i=1}^4  \left| \int y_{h}^{(i)} d\tilde{\rho}_N \right| 
\\ \leq \frac{2C}{N^{\frac{3}{4}}} \left[\sum_{h \geq 2} \frac{1}{\left(1+h^2\right)^r} \sum_{i=1}^4  \left( \int y_{h}^{(i)} d\tilde{\rho}_N \right)^2 \right]^{\frac{1}{2}}
\end{multline*}
with $C^2 := \sum_{h} \frac{h^2}{\left(1+h^2\right)^r}$, where we have used the Cauchy-Schwartz inequality. Summing up, for $t \leq \tau_N^M$,
\[
 \left(\frac{\partial}{\partial x_j} \left\| \tilde{\rho}_N\right\|^2_r\right) ^2 \leq \frac{C^2}{N^{\frac{3}{2}}}\left\| \tilde{\rho}_N\right\|^2_r \leq \frac{C^2 M}{N^{\frac{3}{2}}},
 \]
from which \eqref{Kcollaps2} follows.
\end{proof}

\begin{remark}
For later use, we observe that the $M$-dependence of the constant $C$ in \eqref{Kcollaps1} comes form the estimates in \eqref{Kkeypoint} and \eqref{Kcollaps7}, where the factor $Y_1^{(j,N)} $ is estimated by a constant $C(M)$. If we replace such estimate with the trivial one
\[
\left|Y_1^{(j,N)} \right| \leq N^{\frac{1}{4}},
\]
we obtain the following estimate, which does not require any stopping argument:
\[
J_N \left\| \tilde{\rho}_N \right\|^2_r (\sqrt{N} t) \leq - \left(\frac{1}{2} - 2 \omega^2\right) \sqrt{N} \left\| \tilde{\rho}_N \right\|^2_r (\sqrt{N} t) + C \sqrt{N},
\]
which implies
\begin{equation}\label{Kkeypoint2}
\sup_{N \geq 1, t \geq 0} E\left[ \left\| \tilde{\rho}_N \right\|^2_r \right] < +\infty.
\end{equation}
\end{remark}

We can now prove the main result of this section, corresponding to the first part of Theorem \ref{Kinhom}.
\begin{proposition} \label{Kcollapsfull}
For every $T>0$, the sequence $\left(\left\| \tilde{\rho}_N\right\|_r (\sqrt{N} t) \right)_{t \in [0,T]}$ collapses to zero.
\end{proposition}
\begin{proof}
Given the result of Lemma \ref{12Kinh}, all we have to show is that for every $\varepsilon >0$ there exist $M,N_0 >0$ such that
\begin{equation} \label{Kcollaps13}
\sup_{N \geq N_0}P \{ \tau_N^M > T \} \leq \varepsilon.
\end{equation}
Consider the function 
\[
\psi(x,y) := \sqrt{1 + x^2 + y^2}.
\]
Note that $\psi$ has uniformly bounded partial derivatives, and $\psi(x,y) \geq \min(|x|,|y|)$. We begin by observing that
\begin{multline} \label{Kcollaps14}
\left\{\tau_N^M > T \right\} \subseteq \left\{\sup_{0 \leq t \leq T \wedge \tau_N^M} \left\| \tilde{\rho}_N\right\|^2_r (\sqrt{N} t) \geq \frac{M}{2} \right\} \\
\cup  \left\{\sup_{0 \leq t \leq T \wedge \tau_N^M} \psi\left(V_1^{(1,N)}(t), V_1^{(2,N)}(t)\right) \geq \frac{M}{2} \right\}.
\end{multline}
By Lemma \ref{12Kinh}, for $N$ large the probability 
\[
P\left\{ \sup_{0 \leq t \leq T \wedge \tau_N^M} \left\| \tilde{\rho}_N\right\|^2_r (\sqrt{N} t) \geq \frac{M}{2} \right\}
\]
 can be made arbitrarily small. Thus, \eqref{Kcollaps13} follows if we show that for every $\varepsilon >0$ there exist $M,N_0 >0$ such that
 \begin{equation} \label{Kcollaps15}
\sup_{N \geq N_0}P\left\{ \sup_{0 \leq t \leq T \wedge \tau_N^M} \psi\left(V_1^{(1,N)}(t), V_1^{(2,N)}(t)\right) \geq \frac{M}{2} \right\} \leq \varepsilon \,.
\end{equation}
For the proof of \eqref{Kcollaps15} we consider the perturbation 
\[
\psi_N = \psi + N^{-\frac{1}{4}} \psi_1
\]
as illustrated in Section \ref{ssK:PerturbationTheory}, with
\[
 \psi_1 := - \left(L^{(1)}\right)^{-1} L^{(2)} \psi.
 \]
As seen in Section \ref{ssK:PerturbationTheory}, $\psi_1$  is a linear combination of terms of the form
\[
\partial_i \psi \left(\int v_1^{(1)} d\tilde{\rho}_N, \int v_1^{(2)} d\tilde{\rho}_N \right) \int v_1^{(j)} d\tilde{\rho}_N  \int v_2^{(h)} d\tilde{\rho}_N ,
\]
and therefore, up to time $\tau_N^M$, can be bounded in absolute value by some $M$-dependent constant $C(M)$. This implies that, for every given $M$ and for large enough $N$
\begin{equation} \label{Kcollaps16}
 P\left\{ \sup_{0 \leq t \leq T \wedge \tau_N^M} \psi\left(V_1^{(1,N)}(t), V_1^{(2,N)}(t)\right) \geq \frac{M}{2} \right\} \leq P\left\{ \sup_{0 \leq t \leq T \wedge \tau_N^M} \psi_N\left(\cdot \right) \geq \frac{M}{3} \right\}.
\end{equation}
By abuse of notation, we write $\psi_N(t)$ in place of 
\[
\psi_N\left(V_h^{(i,N)}(t): h=1,2; \, i=1,2,3,4\right).
\]
 Consider the semimartingale representation
 \begin{equation} \label{Kcollaps17}
\psi_N(t) = \psi_N(0) +  \int_0^t J_N \psi_N (s) ds + \mathcal{M}_N(t),
\end{equation}
where
 \begin{equation} \label{Kcollaps18}
\mathcal{M}_N(t) =  N^{\frac{1}{4}} \sum_{j=1}^N \int_0^t  \frac{\partial \psi_N}{\partial x_j} (s) dW_j(s)
\end{equation}
where $\{W_j(t) : t>0, j=1,\dots, N\}$ is a system of independent standard Brownian motions on $[0,2\pi]$. We have
\begin{multline*}
P\left\{ \sup_{0 \leq t \leq T \wedge \tau_N^M} \psi_N\left(\cdot \right) \geq \frac{M}{3} \right\} \leq P\left\{ \psi_N(0) \geq \frac{M}{9} \right\} \\ 
+ P\left\{ \sup_{0 \leq t \leq T \wedge \tau_N^M} J_N \psi_N (t) \geq \frac{M}{9T} \right\}  + P\left\{ \sup_{0 \leq t \leq T \wedge \tau_N^M} \mathcal{M}_N(t) \geq \frac{M}{9} \right\}.
\end{multline*}
The term $P\left\{\psi_N(0) \geq \frac{M}{9} \right\}$ is easy to control, since the random variables $V_h^{(i,N)}(0)$ converge to zero in probability. We are therefore left to show that the probabilities
 \begin{equation} \label{Kcollaps19}
P\left\{ \sup_{0 \leq t \leq T \wedge \tau_N^M} J_N \psi_N (t) \geq \frac{M}{9T} \right\}
\end{equation}
and
 \begin{equation} \label{Kcollaps20}
P\left\{ \sup_{0 \leq t \leq T \wedge \tau_N^M} \mathcal{M}_N(t) \geq \frac{M}{9} \right\}
\end{equation}
are small for $N$ large enough. We begin to deal with \eqref{Kcollaps19}. By \eqref{Kexpansion} and the choice of $\psi_1$, we have
\[
J_N \psi_N = L^{(3)} \psi + L^{(2)} \psi_1 + o(1),
\]
where the term $o(1)$ is bounded by $\frac{C(M)}{N^{\alpha}}$ for some $\alpha>0$. Moreover, it is easily shown that $L^{(3)} \psi$ is bounded uniformly in $N$ and $M$. To deal with $L^{(2)} \psi_1$, we use Proposition \ref{Kpropperthth}, which gives
\begin{multline}\label{Kcollaps21}
L^{(2)} \psi_1(t) \\ = - \frac{(1+4\omega^2)^2(1-8\omega^2)}{4(1-4\omega^2)^3(1+\omega^2)} V_1^{(1,N)}(t) \left[ \left( V_1^{(1,N)}(t) \right)^2 + \left( V_1^{(2,N)}(t) \right)^2 \right] \partial_1 \psi(\cdot, \cdot) \\ -  \frac{(1+4\omega^2)^2(1-8\omega^2)}{4(1-4\omega^2)^3(1+\omega^2)} V_1^{(2,N)}(t) \left[ \left( V_1^{(1,N)}(t) \right)^2 + \left( V_1^{(2,N)}(t) \right)^2 \right] \partial_2 \psi(\cdot, \cdot) \\ + \mbox{ collapsing terms},
\end{multline}
where, again, the ``collapsing terms'' are bounded by $\frac{C(M)}{N^{\alpha}}$. Observing that 
\[
\partial_i \psi(\cdot, \cdot) = \frac{V^{(i,N)}}{\psi(\cdot, \cdot)},
\]
and since, by assumption, $1-8\omega^2 \geq 0$, the non-collapsing part of \eqref{Kcollaps21} is nonnegative. We therefore conclude that, for $t \leq \tau_N^M$
\[
J_N \psi_N \leq C + \frac{C(M)}{N^{\alpha}}
\]
with $C$ independent of $M,N$. This implies that the probability in \eqref{Kcollaps19} is arbitrarily small for $M$ (first) and $N$ (then) sufficiently large.

\noindent
We now deal with \eqref{Kcollaps20}. By Doob's Maximal Inequality
\begin{align} \label{Kcollaps22}
P\left\{ \sup_{0 \leq t \leq T \wedge \tau_N^M} \mathcal{M}_N(t) \geq \frac{M}{9} \right\} &\leq \frac{E\left[ \left( \mathcal{M}_N(T \wedge \tau_N^M)\right)^2 \right]}{(M/9)^2}  \nonumber\\ 
&= \frac{N^{\frac{1}{2}} \sum_{j=1}^N \int_0^{T \wedge \tau_N^M}  \left(\frac{\partial \psi_N}{\partial x_j} (t) \right)^2dt}{(M/9)^2}.
\end{align}
Up to term bounded by $\frac{C(M)}{N^{\alpha}}$, we can replace $\frac{\partial \psi_N}{\partial x_j}$ with $\frac{\partial \psi}{\partial x_j}$ in \eqref{Kcollaps22}. Moreover
\begin{equation} \label{Kcollaps22bis}
\left|\frac{\partial \psi}{\partial x_j}\right| = \left| \sum_{i=1,2} \partial_i \psi \, \frac{1}{N^{\frac{3}{4}}} \frac{\partial v_1^{(i)}}{\partial x}(x_j,\eta_j) \right| \leq \frac{C}{N^{\frac{3}{4}}} 
\end{equation}
for a constant $C$ independent of $M,N$. Inserting this in \eqref{Kcollaps22}, we have, for some $C,C(M)>0$,
\[
P\left\{ \sup_{0 \leq t \leq T \wedge \tau_N^M} \mathcal{M}_N(t) \geq \frac{M}{9} \right\} \leq \frac{C + \frac{C(M)}{N^{\alpha}}}{(M/9)^2},
\]
which, again, is small for $M$ (first) and $N$ (then) sufficiently large. This completes the proof.
\end{proof}

\begin{remark} \label{Kexpl}
The assumption $\omega \leq \frac{1}{2 \sqrt{2}}$ has been used in \eqref{Kcollaps21}, to obtain bounds for $J_N \psi_N $. When the processes are stopped, as in the part of Theorem \ref{Kinhom} concerning the case $\frac{1}{2 \sqrt{2}} < \omega < \frac{1}{2}$, those estimates are essentially trivial because of the uniform boundedness of the stopped processes.
\end{remark}

\subsection{Identification of the Limiting Generator and Convergence}\label{ssK:Convergence}


In this Section we complete the proof of Theorem \ref{Kinhom}. The argument follows that of Section \ref{ssCW:Convergence}, so most details are omitted.

\noindent
The candidate for the limiting generator in \eqref{Klimgen} has been obtained in Proposition \ref{Kpropperthth} for the {\em drift} part, while the diffusion part comes from the term $L^{(3)} \psi$ that, by direct computation, is shown to be equal to
\[
L^{(3)} \psi = \frac{1+\omega^2}{4} \left[ \partial^2_{11} \psi +  \partial^2_{22} \psi  \right].
\]
In what follows we denote by $J$ the generator of the diffusion process in Theorem \ref{Kinhom}. The proof of convergence develops along the following steps.

\noindent
{\em Step 1: tightness of the processes $V_1^{(1,N)}$ and $V_1^{(2,N)}$}. We use conditions \eqref{T1} and \eqref{T2}. Due to \eqref{Kcollaps13}, we are allowed to stop the processes at $\tau_N^M$ for some large $M$. Condition \eqref{T1} can be obtained simultaneously for $V_1^{(1,N)}$ and $V_1^{(2,N)}$ by \eqref{Kcollaps15}. In order to establish \eqref{T2} for, e.g., $V_1^{(1,N)}$, we consider the function
\[
\psi \left( V_1^{(1,N)},V_1^{(2,N)} \right) := V_1^{(1,N)},
\]
together with its perturbation $\psi_N$ as in \eqref{Kperturb} and \eqref{psi1}. Up to $o(1)$ terms, for stopping times $\tau_1 \leq \tau_2$,
\[
V_1^{(1,N)}(\tau_2) - V_1^{(1,N)}(\tau_1)  \simeq \int_{\tau_1}^{\tau_2} J_N \psi_N dt +  N^{\frac{1}{4}}\int_{\tau_1}^{\tau_2} \sum_{j=1}^N \frac{\partial \psi_N}{\partial x_j} dW_j(t).
\]
As in the proof of Proposition \ref{Kcollapsfull}, we find a (possibly $M$-dependent) constant $C$ such that the uniform bound
\[
\left|J_N \psi_N \right| + N^{\frac{1}{2}}\sum_{j=1}^N \left(\frac{\partial \psi_N}{\partial x_j}\right)^2 \leq C
\]
holds. This implies
\[
\sup_{\tau_1 \leq \tau_2 \leq \tau_1 + \delta} E\left[ \left|V_1^{(1,N)}(\tau_2) - V_1^{(1,N)}(\tau_1) \right| \right]  \leq C \delta
\]
that, by Chebischev inequality, yields \eqref{T2} for $V_1^{(1,N)}$.

\noindent
{\em Step 2: convergence to the solution of a martingale problem}. Denote by $\left( V_1^{(1,n)},V_1^{(2,n)} \right)$ a convergent subsequence of $\left( V_1^{(1,N)},V_1^{(2,N)} \right)$. For a function $\psi: \mathbb{R}^2 \rightarrow \mathbb{R}$ of class ${\cal{C}}^2$ and with bounded derivatives, denote by $\psi_n$ its perturbation  as in \eqref{Kperturb} and \eqref{psi1}. Consider the martingale
\begin{equation} \label{Kmart1}
\mathcal{M}_n(t) := \psi_n(t)- \psi_n(0) - \int_0^t J_n \psi_n(s)ds = n^{\frac{1}{4}} \int_0^t \sum_{j=1}^n \frac{\partial}{\partial x_j} \psi_n(s) dW_j(s).
\end{equation}
It should be recalled that $\psi_n$ is a function of $V_1^{(i,n)},V_2^{(i,n)}$, $i=1,2,3,4$, so when we write $\psi_n(t)$ we mean that $t$ is the time at which the processes in the argument of $\psi_n$ are evaluated. Considering that:
\begin{itemize}
\item
$\psi_n \rightarrow \psi$ as $n \rightarrow +\infty$ uniformly on compact sets;
\item
the processes $V_1^{(i,n)},V_2^{(i,n)}$ admit a weak limit $V_1^{(i)},V_2^{(i)}$, which is zero for $V_1^{(i,n)}$, $i=3,4$ and $V_2^{(i,n)}$, $i=1,2,3,4$,
\end{itemize}
it follows that the process $\mathcal{M}_n(t)$ converges weakly to 
\[
\mathcal{M}(t) = \psi \left( V_1^{(1)}(t), V_1^{(2)}(t) \right) - \psi \left( V_1^{(1)}(0), V_1^{(2)}(0) \right) - \int_0^t J\psi \left( V_1^{(1)}(s), V_1^{(2)}(s) \right)ds.
\]
If we show that, for each $\psi$ with the properties specified above, $\mathcal{M}(t)$ is a martingale, then we have that the limiting processes $\left( V_1^{(1)}(t), V_1^{(2)}(t) \right)$ solve the martingale problem for $J$; since uniqueness holds for this martingale problem, the proof of Theorem \ref{Kinhom} would be completed. It is therefore enough to show that $\mathcal{M}(t)$ is a martingale. Similarly to what we have done in Lemma \ref{31CWG}, it suffices to show that, for every $t >0$,
\[
\sup_n E\left[\left(\mathcal{M}_n(t)\right)^2\right] < +\infty.
\]
Note that
\[
E\left[\left(\mathcal{M}_n(t)\right)^2\right] = n^{\frac{1}{2}} \sum_{j=1}^n \int_0^t E\left[ \left( \frac{\partial}{\partial x_j} \psi_n(s)\right)^2 \right] ds.
\]
Thus, it is enough to show that, for some constant $C>0$, the inequality
\begin{equation} \label{Kmart2}
E\left[ \left( \frac{\partial}{\partial x_j} \psi_n(s)\right)^2 \right] \leq \frac{C}{n^{\frac{3}{2}}}
\end{equation}
is satisfied. \\
It should be noticed that in \eqref{Kcollaps22bis} we gave a {\em pointwise} estimate (i.e. not in mean) of this sort; that, however, holds for the {\em unperturbed} function $\psi$. In that case the difference between $\psi$ and its perturbation $\psi_n$ was estimated by a bound of the form $\frac{C(M)}{N^{\alpha}}$. But now we are not stopping the process anymore, so a little more care is needed. We recall that
\[
\psi_n = \psi + n^{-\frac{1}{4}} \psi_1.
\]
Given the bound in  \eqref{Kcollaps22bis}, in order to obtain \eqref{Kmart2} it is enough to show that
\begin{equation} \label{Kmart3}
E\left[ \left( \frac{\partial}{\partial x_j} \psi_1(s)\right)^2 \right] \leq \frac{C}{n}.
\end{equation}
As seen in Section \ref{ssK:PerturbationTheory}, $\psi_1$  is a linear combination of terms of the form
\[
F := \partial_i \psi \left(\int v_1^{(1)} d\tilde{\rho}_n, \int v_1^{(2)} d\tilde{\rho}_n \right) \int v_1^{(l)} d\tilde{\rho}_n  \int v_2^{(h)} d\tilde{\rho}_n,
\]
$i=1,2$, $l,h = 1,2,3,4$.
So it is enough to consider one of such terms. We have
\begin{multline*}
\frac{\partial F}{\partial x_j} = \frac{1}{n^{\frac{3}{4}}} \left[ \partial^2_{1,i} \psi \left(\int v_1^{(1)} d\tilde{\rho}_n, \int v_1^{(2)} d\tilde{\rho}_n \right) \frac{\partial}{\partial x}v_1^{(1)} (x_j,\eta_j) \int v_1^{(j)} d\tilde{\rho}_n  \int v_2^{(h)} d\tilde{\rho}_n \right. \\ \left. +  \partial^2_{2,i} \psi \left(\int v_1^{(1)} d\tilde{\rho}_n, \int v_1^{(2)} d\tilde{\rho}_n \right) \frac{\partial}{\partial x}v_1^{(2)} (x_j,\eta_j) \int v_1^{(j)} d\tilde{\rho}_n  \int v_2^{(h)} d\tilde{\rho}_n \right. \\ + \partial_i \psi \left(\int v_1^{(1)} d\tilde{\rho}_n, \int v_1^{(2)} d\tilde{\rho}_n \right) \frac{\partial}{\partial x}v_1^{(l)} (x_j,\eta_j) \int v_2^{(h)} d\tilde{\rho}_n \\ \left. + \partial_i \psi \left(\int v_1^{(1)} d\tilde{\rho}_n, \int v_1^{(2)} d\tilde{\rho}_n \right) \frac{\partial}{\partial x}v_2^{(h)} (x_j,\eta_j)  \int v_1^{(l)} d\tilde{\rho}_n \right].
\end{multline*}
Consider the first of the summands above, the others can be dealt with similarly. The factor
\[
 \partial^2_{1,i} \psi \left(\int v_1^{(1)} d\tilde{\rho}_n, \int v_1^{(2)} d\tilde{\rho}_n \right) \frac{\partial}{\partial x}v_1^{(1)} (x_j,\eta_j)
 \]
is uniformly bounded. Also the term
\[
\frac{1}{n^{\frac{1}{4}}} \int v_1^{(j)} d\tilde{\rho}_n 
\]
is uniformly bounded.
The last factor, $\int v_2^{(h)} d\tilde{\rho}_n $, is clearly bounded in absolute value by $ \left\| \tilde{\rho}_n \right\|_r $, defined in \eqref{Knormr}. Estimating similarly all terms, one sees that
\[
E\left[\left(\frac{\partial F}{\partial x_j} \right)^2\right] \leq \frac{C}{n} E\left[ \left\| \tilde{\rho}_n \right\|^2_r \right] \leq \frac{C'}{n}
\]
for some constants $C,C'$, where we have used \eqref{Kkeypoint2}. This establishes \eqref{Kmart3}, and thus completes the proof of Theorem  \ref{Kinhom}.



\begin{thebibliography}{99}\addcontentsline{toc}{section}{References}
\bibitem{MaPe91} Joao M. G. Amaro de Matos  and J. Fernando Perez. Fluctuations in the {C}urie-{W}eiss version of the random field {I}sing model. \emph{J. Statist. Phys.}, 62:587--608, 1991.
\bibitem{BrDu01} William A. Brock  and Steven N. Durlauf. Discrete choice with social interactions. \emph{Rev. Econom. Stud.}, 68:235--260, 2001.
\bibitem{Col09} Francesca Collet. \emph{The impact of disorder in the critical dynamics of mean-field models.} PhD thesis, Department of Pure and Applied Mathematics, University of Padova, 2009.
\bibitem{CDS10} Francesca Collet, Paolo Dai Pra, Elena Sartori, A Simple Mean Field Model for Social Interactions:
Dynamics, Fluctuations, Criticality, \emph{J. Statist. Phys.}, Volume 139, Number 5, 820-858, 2010.
\bibitem{CoEi88} Francis Comets and Thomas Eisele. Asymptotic dynamics, noncritical and critical fluctuations for a geometric long-range interacting model. \emph{Comm. Math. Phys.}, 118:531--567, 1988.
\bibitem{DaPdHo95} Paolo Dai Pra and Frank den Hollander. Mc{K}ean-{V}lasov limit for interacting random processes in random media. Technical report, Department of Mathematics, University of Nijmegen, 1995.
\bibitem{DPRST09} Paolo Dai Pra, Wolfgang J. Runggaldier, Elena Sartori, and Marco Tolotti. Large portfolio losses; A dynamic contagion model. \emph{Ann. Appl. Probab.}, 19:347--394, 2009.
\bibitem{DaPTo09} Paolo Dai Pra  and Marco Tolotti. Heterogeneous credit portfolios and the dynamics of the aggregate losses. \emph{Stochastic Processes Appl.}, 119:2913--2944, 2009.
\bibitem{Daw83} Donald A. Dawson, Critical dynamics and fluctuations for a mean-field model of cooperative behavior. \emph{J. Statist. Phys.}, Volume 31, Number 1, 29-85, 1983.
\bibitem{DGGaLa} Corrado Di Guilmi, Mauro Gallegati, and Simone Landini. Financial fragility, mean-field interaction and macroeconomic dynamics: a stochastic model. Preprint. Available at SSRN: \textsf{http://ssrn.com/abstract=1258542}, 2008.
\bibitem{EtKu86} Stewart N. Ethier  and  Thomas G. Kurtz. \emph{Markov processes: characterization and convergence}. John Wiley \& Sons Inc., New York, 1986.
\bibitem{For09} Marco Formentin. \emph{Two problems concerning interacting systems: 1. {O}n the {P}urity of the free boundary condition {P}otts measure on {G}alton-{W}atson trees 2. {U}niform propagation of chaos and fluctuation theorems in some spin-flip models.} PhD thesis, Department of Pure and Applied Mathematics, University of Padova, 2009.
\bibitem{FrBa08} R{\"u}diger Frey  and Jochen Backhaus. Pricing and hedging of portfolio credit derivatives with interacting default intensities. \emph{Int. J. Theor. Appl. Finance}, 11:611--634, 2008.
\bibitem{Krt73} Thomas G. Kurtz. A limit theorem for perturbed operator semigroups with applications to random evolutions. \emph{J. Functional Analysis}, 12:55--67, 1973.
\bibitem{LaLi07} Jean-Michel Lasry and Pierre-Louis Lions. Mean field games. \emph{Jpn. J. Math.}, 2:229--260, 2007.
\bibitem{Pap77} George C. Papanicolaou. Introduction to the asymptotic analysis of stochastic equations. In \emph{Lectures in Applied Mathematics 16}, pages 108--148. American Mathematical Society, 1977.
\bibitem{PaStVa77} George C. Papanicolaou, Daniel Stroock, and S. R. S. Varadhan. Martingale approach to some limit theorems. In \emph{Papers from the {D}uke {T}urbulence {C}onference ({D}uke {U}niv., {D}urham, {N}.{C}., 1976), {P}aper {N}o. 6}. Duke University, 1977.
\bibitem{Sar07} Elena Sartori. \emph{Some aspects of spin systems with mean-field interaction}. PhD thesis, Department of Pure and Applied Mathematics, University of Padova, 2007.
\bibitem{Shi96} Albert N. Shiryaev. \emph{Probability}. Springer-Verlag, New York, second edition, 1996.
\end{thebibliography}

\end{document}